\documentclass[preprint,11pt,nopreprintline]{elsarticle}

\usepackage{amsmath,amssymb,amsthm,amsfonts,amsopn}
\usepackage{graphicx}
\usepackage{multirow}
\usepackage{subfigure}
\usepackage{mathrsfs}
\usepackage[margin=1in]{geometry}
\usepackage[hidelinks]{hyperref}
\usepackage{cleveref}
\usepackage{lineno}

\crefname{equation}{equation}{equations}

\journal{Journal of Computational Physics}

\theoremstyle{plain}
\newtheorem{thm}{Theorem}[section]

\theoremstyle{definition}

\theoremstyle{remark}
\newtheorem{remark}[thm]{Remark}

\newcommand{\paren}[1]{\left(#1\right)}

\newcommand{\bm}[1]{\boldsymbol{#1}}

\hypersetup{
  pdftitle={ADI Schemes for Three-Dimensional Heat Equations with Irregular Boundaries and Interfaces},
  pdfauthor={Han Zhou, Minsheng Huang, and Wenjun Ying}
}

\begin{document}
% \linenumbers

\begin{frontmatter}

\title{ADI Schemes for Three-Dimensional Heat Equations with Irregular Boundaries and Interfaces}

\author[inst1]{Han Zhou}
\ead{zhouhan@sjtu.edu.cn}

\author[inst1]{Minsheng Huang}
\ead{mingo.stemon@sjtu.edu.cn}

\author[inst1]{Wenjun Ying\corref{cor1}}
\ead{wying@sjtu.edu.cn}

\cortext[cor1]{Corresponding author}

\address[inst1]{School of Mathematical Sciences, MOE-LSC and Institute of Natural Sciences, Shanghai Jiao Tong University, Minhang, Shanghai 200240, China}

\begin{abstract}
We develop efficient alternating direction implicit (ADI) schemes for three-dimensional heat equations with irregular boundaries and interfaces. Starting from the classical Douglas-Gunn ADI scheme, we construct a modified variant that removes the accuracy degradation commonly observed in the presence of time-dependent boundary conditions. The unconditional stability of the modified scheme is established by Fourier analysis.
By combining the ADI framework with a one-dimensional kernel-free boundary integral (KFBI) method, we obtain KFBI-ADI schemes for heat equations on irregular domains. In the resulting one-dimensional subproblems, the KFBI discretization preserves the structure of the coefficient matrix on Cartesian grids and therefore permits an efficient implementation based on the Thomas algorithm.
Several numerical experiments for the heat equation and a reaction-diffusion equation confirm second-order accuracy and unconditional stability. For the Stefan problem, we couple the ADI framework with a level set method to capture the moving interface. Numerical examples of three-dimensional dendritic solidification are also presented.
\end{abstract}


\end{frontmatter}

\section{Introduction}
\label{sec:intro}
The heat equation and related problems, such as reaction-diffusion equations and free-boundary problems, arise in a broad range of applications, including heat transfer \cite{YANENKO19631094, Gibou2005, Pandit2016}, pattern formation \cite{Chen1992, Fernandes2012, Asante2020}, crystal growth \cite{Soni1999, Nochetto1991, Nochetto1991b}, and many other areas.
The development of efficient numerical algorithms for such problems is therefore of substantial interest, particularly in three spatial dimensions and in the presence of complex boundaries or interfaces.

ADI schemes are well known to be highly efficient for solving different kinds of partial differential equations (PDEs) in multiple spatial dimensions \cite{doi:10.1137/0103004,doi:10.1137/0103003, Douglas1964, Kim2007, Zhao2014}.
Compared with conventional explicit and implicit time-stepping methods, which either suffer from severe stability constraints or require solving large algebraic systems, ADI schemes are generally unconditionally stable and have computational costs comparable to those of explicit schemes due to the dimension-splitting strategy.
Classical ADI schemes are primarily designed for simple geometries that can be discretized by Cartesian grids, such as tensor-product domains (e.g. rectangles in 2D and cubes in 3D) or unions of such domains (e.g. an L-shaped domain).
In more realistic settings, however, domain boundaries and internal interfaces are often geometrically complex, which makes the direct application of Cartesian-grid ADI schemes difficult.
One of the pioneering developments for handling complex boundaries on Cartesian grids is the immersed boundary method (IBM), introduced by C. S. Peskin in the 1970s \cite{Peskin1977, Peskin2002}.
Motivated by the success of IBM, many researchers have studied Cartesian grid methods, and numerous approaches have emerged over the past decades \cite{doi:10.1137/0721021,doi:10.1137/0731054,doi:10.1137/1.9780898717464, ZHOU20061}.
Because of the simplicity and efficiency of dimension splitting, various ADI schemes have been proposed for PDEs with complex interfaces, including immersed interface method (IIM)-based ADI schemes \cite{Li1994, IIMBook, Liu2014, Li2018}, matched interface and boundary (MIB)-based ADI schemes \cite{Zhao2015, Wei2018, Li2021}, ghost fluid method (GFM)-based ADI schemes \cite{Li2020, Li2021} for parabolic interface problems, and ADI Yee schemes for Maxwell's equations with discontinuous coefficients \cite{Deng2021}.
There have also been attempts to solve nonlinear interface problems, such as the IIM-ADI scheme for nonlinear convection-diffusion equations \cite{Liu2013} and the ADI scheme for the nonlinear Poisson-Boltzmann equation \cite{Geng2013}.

Recently, several dimension-splitting methods based on the 1D KFBI method have been proposed for solving time-dependent PDEs on complex domains in two spatial dimensions \cite{Zhou2023}.
The KFBI method was originally proposed by W. Ying and C. S. Henriquez in 2007 for solving elliptic PDEs with irregular boundaries in two spatial dimensions \cite{Ying2007} and was generalized to three spatial dimensions \cite{Ying2013} and later extended to a fourth-order version \cite{Xie2020}.
The KFBI method has also been successfully applied to interface problems with variable coefficients \cite{Ying2014}, incompressible flows \cite{Xie2020a}, and singularly perturbed reaction-diffusion equations \cite{XIE2021298}.

In this work, we extend the KFBI-ADI approach to the heat equation, reaction-diffusion equations, and the Stefan problem with complex, and possibly moving, boundaries or interfaces in three spatial dimensions.
First, based on the Douglas-Gunn (DG) scheme \cite{Douglas1962}, we introduce a modified Douglas-Gunn (mDG) scheme that improves the accuracy of DG-ADI discretizations for problems with time-dependent boundary conditions.
We then combine both ADI schemes with the 1D KFBI method to obtain geometrically flexible solvers for the heat equation and reaction-diffusion equations on complex domains.
The resulting methods achieve second-order accuracy in both space and time, as demonstrated numerically on a range of test problems.
Finally, for the Stefan problem, we use a level set method \cite{Sethian1985, OSHER198812} to capture the moving boundary and solve the associated heat-equation interface problem within the ADI framework.

The remainder of the paper is organized as follows.
In \Cref{sec:problem}, we describe the model problems considered in this work.
Then, the construction of ADI schemes is described in \Cref{sec:method}.
Some numerical examples are presented in \Cref{sec:result}.
Finally, we provide a brief discussion of the proposed method in \Cref{sec:discu}.

\section{Governing equations}\label{sec:problem}
Let $\Omega \subset \mathbb{R}^3$ be a fixed and complex-shaped domain with boundary $\Gamma = \partial \Omega$.
Let $u(\bm{x}, t): \Omega \times [0,T] \to \mathbb{R}$ denote the unknown temperature field.
We consider the Dirichlet boundary value problem (BVP) for the heat equation:
\begin{subequations}
\label{eqn:heat-eqn}
\begin{align}
    \partial_t u(\bm{x}, t) &= \Delta u(\bm{x}, t) + f(\bm{x}, t), & (\bm{x}, t) \in \Omega \times [0, T], \\
    u(\bm{x}, t) &= g_D(\bm{x}, t), & (\bm{x}, t) \in \partial \Omega \times [0, T], \\
    u(\bm{x}, 0) &= u_0(\bm{x}), & \bm{x} \in \Omega,
\end{align}
\end{subequations}
where $f(\bm{x}, t)$ is the source term, and $g_D(\bm{x}, t)$ and $u_0(\bm{x})$ are the prescribed boundary and initial conditions, respectively.

When $u(\bm{x}, t)$ represents other physical quantities, such as the concentration of a chemical species, equation \eqref{eqn:heat-eqn} is often referred to as the diffusion equation.
In many cases, the source term $f(\bm{x}, t, u)$ depends on $u$, representing chemical reactions in the system.
In this case, equation \eqref{eqn:heat-eqn} becomes a reaction-diffusion equation.

An important application of the heat equation is the Stefan problem, which models phase transitions between solid and liquid phases. In this setting, the solid-liquid interface, also called the free boundary, is unknown a priori and must be determined together with the temperature field. The Stefan problem is therefore a prototypical free-boundary problem.

Let $\Gamma(t)$ represent the free boundary separating the solid region $\Omega_s(t)$ and the liquid region $\Omega_l(t)$. The Stefan problem is formulated as follows \cite{Schmidt1996,Chen1997}:
\begin{subequations}
\begin{align}
    \partial_t u(\bm{x}, t) &= \Delta u(\bm{x}, t), &(\bm{x}, t) \in \Omega_s(t) \cup \Omega_l(t) \times [0,T], \\
    u(\bm{x}, t) &= -\varepsilon_C H - \varepsilon_V V_n, &(\bm{x}, t) \in \Gamma(t) \times [0,T], \label{eqn:gibbs} \\
    V_n &= [D \partial_{\bm{n}} u], &(\bm{x}, t) \in \Gamma(t) \times [0,T], \label{eqn:ste-eqn}
\end{align}
\end{subequations}
subject to initial conditions for both the temperature and the free boundary:
\begin{align}
    u(\bm{x}, 0) &= u_0(\bm{x}), \text{ for }\bm{x}\in \Omega_s(0)\cup\Omega_l(0),  \quad
    \Gamma(0) = \Gamma_0, 
\end{align}
where $H$ is the mean curvature of $\Gamma$ (the sum of principal curvatures), $V_n$ is the normal velocity of the free boundary $\Gamma$, and $\varepsilon_C$ and $\varepsilon_V$ are the surface tension and molecular kinetic coefficients, respectively. The unit outward normal to the boundary is denoted by $\bm{n}$.

The notation $[q] = q_s - q_l$ represents the jump in the one-sided limit values of the quantity $q$ across $\Gamma$. Equations \eqref{eqn:gibbs} and \eqref{eqn:ste-eqn} correspond to the Gibbs-Thomson relation and the Stefan equation, respectively, which couple the temperature with the motion of the free boundary.

In anisotropic solidification problems, the coefficients $\varepsilon_C$ and $\varepsilon_V$ depend on the orientation of the free boundary and can be expressed as $\varepsilon_C(\bm{n})$ and $\varepsilon_V(\bm{n})$.

In this work, we consider the Stefan problem within a bounded region $\mathcal{B} = \Omega_s \cup \Gamma \cup \Omega_l$. We assume that $\mathcal{B}$ is a cubic domain whose boundaries are sufficiently far from $\Gamma$, so that $\partial\mathcal{B} = \partial\Omega_l \setminus \Gamma$. The boundary condition for $u$ on $\partial\mathcal{B}$ is taken to be the homogeneous Neumann condition $\partial_{\bm{n}} u = 0$.

\section{Numerical method}\label{sec:method}
\subsection{ADI schemes}
In this subsection, we describe the alternating direction implicit method used for the time discretization of the heat equation \eqref{eqn:heat-eqn}.
The time interval $[0,T]$ is uniformly partitioned into $N_t$ intervals $0=t^0<t^1\cdots < t^{N_t-1} < t^{N_t} = T$ with constant time step $\tau = T/N_t$. 
The semi-discrete approximation to the exact solution $u(\bm{x}, t^n)$ without spatial discretization is denoted as $u^n$.
\subsubsection{The Douglas-Gunn scheme}
The classical second-order DG scheme \cite{Douglas1962,Hundsdorfer2003} for the 3D heat equation \eqref{eqn:heat-eqn} without source or reaction terms (i.e., $f(\bm{x},t, u) = 0$) is given by
\begin{subequations}
\label{eqn:DG-adi}
    \begin{align}
    \dfrac{u^{\star} - u^n}{\tau} &= \dfrac{1}{2}(u^{\star}_{xx} + u^n_{xx}) + u_{yy}^n + u_{zz}^n, \\
    \dfrac{u^{\star\star} - u^{\star}}{\tau} &= \dfrac{1}{2}(u_{yy}^{\star\star} - u_{yy}^n),\\
    \dfrac{u^{n+1} - u^{\star\star}}{\tau} &= \dfrac{1}{2}(u_{zz}^{n+1} - u_{zz}^n),
    \end{align}
\end{subequations}
where $u^{\star}$ and $u^{\star\star}$ are intermediate variables. 
The DG scheme can be extended directly to problems with $s \geq 3$ operators, so nonzero sources (i.e., $f(t,\bm{x}, u) \neq 0$) can be treated in the same manner \cite{Hundsdorfer2003}.
Theoretical analysis of the DG scheme shows that it is unconditionally stable so that no stability constraint on time steps, such as $\tau<C h^2$ for explicit schemes, is required \cite{Hundsdorfer2003}.

The DG scheme is built as a high-order perturbation of the second-order Crank-Nicolson scheme
\begin{equation}
    \dfrac{u^{n+1} - u^n}{\tau} = \dfrac{1}{2} (\Delta u^{n+1} + \Delta u^n).
\end{equation}
In principle, the intermediate variables $u^{\star}$ and $u^{\star\star}$ are only auxiliary variables, which have little physical meaning.
If one regards $u^{\star}$ and $u^{\star\star}$ as numerical approximations of the exact solution $u(\bm{x}, t^{n+1})$, the scheme is actually a predictor-corrector scheme, which can be shown by rewriting \eqref{eqn:DG-adi} as
\begin{subequations}
\label{eqn:DGADI-equiv}
\begin{align}
\dfrac{u^{\star} - u^{n}}{\tau} &=  \dfrac{1}{2}(u^{\star}_{xx} + u^n_{xx}) + u_{yy}^n + u_{zz}^n,\\
\dfrac{u^{\star\star} - u^n}{\tau} &=  \dfrac{1}{2}(u^{\star}_{xx} + u^n_{xx}) + \dfrac{1}{2}(u_{yy}^{\star\star}+u_{yy}^n) + u_{zz}^n,\\
\dfrac{u^{n+1} - u^{n}}{\tau} &= \dfrac{1}{2}(u^{\star}_{xx} + u^n_{xx}) + \dfrac{1}{2}(u_{yy}^{\star\star}+u_{yy}^n) + \dfrac{1}{2}(u_{zz}^{n+1} + u_{zz}^n).
\end{align}
\end{subequations}
Let $\mathcal{A}_s(u^n, u^{n+1},  u^{\star}, u^{\star\star})$ for $s=1,2,3$ denote difference operators in the three sub-steps. For example,
\begin{equation}
    \mathcal{A}_1(u^n, u^{n+1},  u^{\star}, u^{\star\star}) = \dfrac{u^{\star} - u^{n}}{\tau} - (  \dfrac{1}{2}(u^{\star}_{xx} + u^n_{xx}) + u_{yy}^n + u_{zz}^n).
\end{equation}
Then the DG scheme can be written as
\begin{subequations}
\begin{align}
    \mathcal{A}_1(u^n, u^{n+1},  u^{\star}, u^{\star\star}) &= 0, \\
    \mathcal{A}_2(u^n, u^{n+1},  u^{\star}, u^{\star\star}) &= 0, \\
    \mathcal{A}_3(u^n, u^{n+1},  u^{\star}, u^{\star\star}) &= 0.
\end{align}
\end{subequations}
By replacing $u^n$, $u^{n+1}$,  $u^{\star}$, and $u^{\star\star}$ with exact solutions that they approximate, it can be shown that local truncation errors satisfy
\begin{subequations}\label{eqn:DGADI-LTE}
\begin{align}
E_1^{n+\frac{1}{2}} = \mathcal{A}_1(u(\bm{x},t^n), u(\bm{x},t^{n+1}),  u(\bm{x},t^{n+1}), u(\bm{x},t^{n+1})) &= \mathcal{O}(\tau),\\
E_2^{n+\frac{1}{2}} = \mathcal{A}_2(u(\bm{x},t^n), u(\bm{x},t^{n+1}),  u(\bm{x},t^{n+1}), u(\bm{x},t^{n+1})) &= \mathcal{O}(\tau),\\
E_3^{n+\frac{1}{2}} = \mathcal{A}_3(u(\bm{x},t^n), u(\bm{x},t^{n+1}),  u(\bm{x},t^{n+1}), u(\bm{x},t^{n+1})) &= \mathcal{O}(\tau^2).
\end{align}
\end{subequations}
The insufficient consistency of the first two sub-steps is often overlooked, yet it may reduce the numerical accuracy and convergence order when the boundary conditions for $u^{\star}$ and $u^{\star\star}$ are prescribed using $u^{n+1}$.
This phenomenon will be demonstrated through a numerical example at the end of this subsection.

In fact, there is a natural way to provide boundary conditions for $u^{\star}$ and $u^{\star\star}$ when the domain boundary is axis-aligned, for example, when $\Omega=(0, 1)^3$.
By reversing the order of \eqref{eqn:DG-adi} and restricting the solution to the boundary $\partial\Omega$, boundary conditions for $u^{\star}$ and $u^{\star\star}$ can be computed with a back-tracing approach, also known as the boundary correction technique:
\begin{subequations}\label{eqn:bdry-crc}
\begin{align}
    u^{\star\star} &= (1-\dfrac{\tau}{2}\partial_{zz}) u^{n+1} + \dfrac{\tau}{2}\partial_{zz} u^n, \\
    u^{\star} &= (1-\dfrac{\tau}{2}\partial_{yy})u^{\star\star} + \dfrac{\tau}{2}\partial_{yy}u^n \nonumber\\
    & = (1-\dfrac{\tau}{2}\partial_{yy})(1-\dfrac{\tau}{2}\partial_{zz}) u^{n+1} + \dfrac{\tau}{2}((1-\dfrac{\tau}{2}\partial_{yy})\partial_{zz} + \partial_{yy}) u^n. 
\end{align}
\end{subequations}
Since the boundary is axis-aligned, the boundary condition for $u^{\star\star}$ on $\partial\Omega\cap\{y=0, 1\}$ and that for $u^{\star}$ on $\partial\Omega\cap\{x=0, 1\}$ can be exactly or approximately computed only using boundary conditions of $u(\bm{x}, t)$ at $t^n$ and $t^{n+1}$.
The boundary correction technique is able to retain the full accuracy of the DG scheme.
Unfortunately, the technique is limited to axis-aligned boundaries and there is no straightforward extension for general irregular domains.

\subsubsection{A modified Douglas-Gunn scheme}
To address the loss of accuracy in the original DG scheme, we modify the scheme so that it is better suited to problems with irregular boundaries.
Since the loss of accuracy results from the insufficient consistency of the first two sub-steps, we modify these sub-steps in \eqref{eqn:DGADI-equiv} using an extrapolation technique:
\begin{subequations}
\label{eqn:mDGADI-equiv}
\begin{align}
\dfrac{u^{\star} - u^n}{\tau} &= \dfrac{1}{2}(u^{\star}_{xx} + u^n_{xx}) + \dfrac{3}{2}u_{yy}^n-\dfrac{1}{2}u_{yy}^{n-1} + \dfrac{3}{2}u_{zz}^n-\dfrac{1}{2}u_{zz}^{n-1}, \\
\dfrac{u^{\star\star} - u^n}{\tau} &= \dfrac{1}{2}(u^{\star}_{xx} + u^n_{xx}) + \dfrac{1}{2}(u_{yy}^{\star\star}+u_{yy}^{n}) + \dfrac{3}{2}u_{zz}^n-\dfrac{1}{2}u_{zz}^{n-1}, \\
\dfrac{u^{n+1} - u^n}{\tau} &= \dfrac{1}{2}(u^{\star}_{xx} + u^n_{xx}) + \dfrac{1}{2}(u_{yy}^{\star\star}+u_{yy}^{n}) + \dfrac{1}{2}(u_{zz}^{n+1}+u_{zz}^{n}).
\end{align}
\end{subequations}
The first two sub-steps are both second-order consistent.
After rearranging terms, the scheme can be written as
\begin{subequations}
\label{eqn:3L-ADI}
    \begin{align}
    \dfrac{u^{\star} - u^n}{\tau} &= \dfrac{1}{2}(u^{\star}_{xx} + u^n_{xx}) + \dfrac{3}{2}u_{yy}^n-\dfrac{1}{2}u_{yy}^{n-1} + \dfrac{3}{2}u_{zz}^n-\dfrac{1}{2}u_{zz}^{n-1}, \\
    \dfrac{u^{\star\star} - u^{\star}}{\tau} &= \dfrac{1}{2}u_{yy}^{\star\star}+\dfrac{1}{2}u_{yy}^{n-1} - u_{yy}^n,\\
     \dfrac{u^{n+1} - u^{\star\star}}{\tau} &= \dfrac{1}{2}u_{zz}^{n+1}+\dfrac{1}{2}u_{zz}^{n-1} - u_{zz}^n.
    \end{align}
\end{subequations}
The modified scheme is a predictor-corrector scheme, and each sub-step is a second-order consistent approximation to the continuous problem.
In the modified scheme, the physical boundary conditions at $t^{n+1}$ can be imposed on $u^{\star}$ and $u^{\star\star}$ without sacrificing second-order accuracy, which makes the method more suitable for irregular domains where the boundary correction technique is not available.
As with the original DG scheme, the modified DG scheme is also unconditionally stable, so no restrictive high-order time-step constraint is required; see \Cref{sec:analysis}.
For a non-vanishing source term $f(\bm{x}, t, u)$, the scheme can also be generalized to
\begin{subequations}
\label{eqn:3L-ADI2}
    \begin{align}
    \dfrac{u^{\star} - u^n}{\tau} &= \dfrac{1}{2}(f^{\star} + f^n) + \dfrac{3}{2}u^n_{xx} - \dfrac{1}{2}u_{xx}^{n-1} \\
    & + \dfrac{3}{2}u_{yy}^n-\dfrac{1}{2}u_{yy}^{n-1} + \dfrac{3}{2}u_{zz}^n-\dfrac{1}{2}u_{zz}^{n-1}, \nonumber\\
    \dfrac{u^{\star\star} - u^{\star}}{\tau} &= \dfrac{1}{2}u_{xx}^{\star\star}+\dfrac{1}{2}u_{xx}^{n-1} - u_{xx}^n,\\
    \dfrac{u^{\star\star\star} - u^{\star\star}}{\tau} &= \dfrac{1}{2}u_{yy}^{\star\star\star}+\dfrac{1}{2}u_{yy}^{n-1} - u_{yy}^n,\\
     \dfrac{u^{n+1} - u^{\star\star\star}}{\tau} &= \dfrac{1}{2}u_{zz}^{n+1}+\dfrac{1}{2}u_{zz}^{n-1} - u_{zz}^n.
    \end{align}
\end{subequations}
where $f^{\star}=f(\bm{x}, t^{n+1}, u^\star)$ and $f^{n}=f(\bm{x}, t^{n}, u^n)$.

The sub-steps of the ADI schemes involve explicit terms containing second derivatives of $u^n$ and $u^{n-1}$, which appear on the right-hand side of the one-dimensional subproblems.
When the domain boundary is axis-aligned, these terms can be easily calculated by applying central differences to the grid data of $u^n$ and $u^{n-1}$.
If the domain is irregular, some stencil points of the central difference scheme may lie outside the domain, and one may need to use one-sided interpolation to evaluate these derivatives.
An alternative approach is to solve for $u^{n}$ first and then store the grid values of $u^n_{zz}$ using the relation $u^{n}_{zz} = \frac{2}{\tau}u^{n} + \widetilde{F}$, where $\widetilde{F}$ is the right-hand side of the one-dimensional subproblem.
By rearranging the order of sub-steps, for example $2\rightarrow 3\rightarrow 1$ and $3\rightarrow 1\rightarrow 2$, the values $u^{n}_{xx}$ and $u^n_{yy}$ can also be computed.
In this way, one can easily obtain all the second derivatives of $u^n$ and $u^{n-1}$ without using one-sided interpolation.
However, this approach is more time-consuming, since it requires roughly twice the computational cost due to the extra time-stepping procedures with different sub-step orders.

\subsubsection{Numerical validation}
To demonstrate the effect of inappropriate boundary conditions for intermediate variables on the numerical accuracy and convergence order, we construct a simple example of the 3D heat equation with a time-dependent boundary condition.
The initial condition and boundary condition are chosen such that the exact solution is given by
\begin{equation}
    u(x, y, z, t) = e^{-6t}\sin(\sqrt{2}x)\cos(\sqrt{2}y)\sin(\sqrt{2}z).
\end{equation}
The computational domain is assumed to be a cube $\Omega=(-1.2,1.2)^3$.
Numerical errors are estimated at the final time $T=0.5$.
The problem is solved using the DG scheme both with and without the boundary correction technique, as well as the modified DG scheme.
Spatial derivatives are approximated using the standard three-point central difference scheme.
The errors and convergence orders obtained by the three methods are summarized in \Cref{tab:ADI-cube}.
\begin{table}[htbp]
\centering
\caption{Numerical results of the DG scheme (with and without correction) and the modified DG scheme for the heat equation.}
\label{tab:ADI-cube}
\resizebox{\textwidth}{!}{
\begin{tabular}{|c|c|cccc|cccc|cccc|}
\hline
\multirow{3}{*}{N} &
  \multirow{3}{*}{$\tau$} &
  \multicolumn{4}{c|}{DG without correction} &
  \multicolumn{4}{c|}{DG with correction} &
  \multicolumn{4}{c|}{modified DG} \\ \cline{3-14} 
 &
   &
  \multicolumn{2}{c|}{$L_2$} &
  \multicolumn{2}{c|}{$L_{\infty}$} &
  \multicolumn{2}{c|}{$L_2$} &
  \multicolumn{2}{c|}{$L_{\infty}$} &
  \multicolumn{2}{c|}{$L_2$} &
  \multicolumn{2}{c|}{$L_{\infty}$} \\ \cline{3-14} 
 &
   &
  \multicolumn{1}{c|}{Error} &
  \multicolumn{1}{c|}{Rate} &
  \multicolumn{1}{c|}{Error} &
  Rate &
  \multicolumn{1}{c|}{Error} &
  \multicolumn{1}{c|}{Rate} &
  \multicolumn{1}{c|}{Error} &
  Rate &
  \multicolumn{1}{c|}{Error} &
  \multicolumn{1}{c|}{Rate} &
  \multicolumn{1}{c|}{Error} &
  Rate \\ \hline
16 &
  1/32 &
  \multicolumn{1}{c|}{5.26e-5} &
  \multicolumn{1}{c|}{} &
  \multicolumn{1}{c|}{5.65e-4} &
   &
  \multicolumn{1}{c|}{1.70e-6} &
  \multicolumn{1}{c|}{} &
  \multicolumn{1}{c|}{4.98e-6} &
   &
  \multicolumn{1}{c|}{1.65e-5} &
  \multicolumn{1}{c|}{} &
  \multicolumn{1}{c|}{1.28e-4} &
   \\ \hline
32 &
  1/64 &
  \multicolumn{1}{c|}{1.37e-5} &
  \multicolumn{1}{c|}{1.94} &
  \multicolumn{1}{c|}{2.47e-4} &
  1.19 &
  \multicolumn{1}{c|}{4.32e-7} &
  \multicolumn{1}{c|}{1.98} &
  \multicolumn{1}{c|}{1.21e-6} &
  2.04 &
  \multicolumn{1}{c|}{2.75e-6} &
  \multicolumn{1}{c|}{2.58} &
  \multicolumn{1}{c|}{2.03e-5} &
  2.66 \\ \hline
64 &
  1/128 &
  \multicolumn{1}{c|}{3.51e-6} &
  \multicolumn{1}{c|}{1.96} &
  \multicolumn{1}{c|}{1.15e-4} &
  1.10 &
  \multicolumn{1}{c|}{1.09e-7} &
  \multicolumn{1}{c|}{1.99} &
  \multicolumn{1}{c|}{2.99e-7} &
  2.02 &
  \multicolumn{1}{c|}{5.26e-7} &
  \multicolumn{1}{c|}{2.39} &
  \multicolumn{1}{c|}{4.24e-6} &
  2.26 \\ \hline
128 &
  1/256 &
  \multicolumn{1}{c|}{8.84e-7} &
  \multicolumn{1}{c|}{1.99} &
  \multicolumn{1}{c|}{5.59e-5} &
  1.04 &
  \multicolumn{1}{c|}{2.74e-8} &
  \multicolumn{1}{c|}{1.99} &
  \multicolumn{1}{c|}{7.42e-8} &
  2.01 &
  \multicolumn{1}{c|}{1.13e-7} &
  \multicolumn{1}{c|}{2.22} &
  \multicolumn{1}{c|}{9.65e-7} &
  2.14 \\ \hline
256 &
  1/256 &
  \multicolumn{1}{c|}{2.22e-7} &
  \multicolumn{1}{c|}{1.99} &
  \multicolumn{1}{c|}{2.76e-5} &
  1.02 &
  \multicolumn{1}{c|}{6.87e-9} &
  \multicolumn{1}{c|}{2.00} &
  \multicolumn{1}{c|}{1.85e-8} &
  2.00 &
  \multicolumn{1}{c|}{2.61e-8} &
  \multicolumn{1}{c|}{2.11} &
  \multicolumn{1}{c|}{2.32e-7} &
  2.06 \\ \hline
\end{tabular}
}
\end{table}
The uncorrected DG scheme yields the least accurate results.
Compared with the corrected scheme, the uncorrected DG scheme also exhibits an order reduction in the $L_{\infty}$ norm from second order to first order.
The modified DG scheme provides more accurate results than the uncorrected DG scheme and obtains second-order convergence in both norms. 

\subsection{KFBI-ADI method for irregular boundaries}
Consider an irregular domain $\Omega\subset\mathbb{R}^3$.
To solve PDEs on $\Omega$, we first embed it into a sufficiently large cuboid bounding box $\mathcal{B}$, which works as the computational domain.
For simplicity, the box $\mathcal{B}$ is assumed to be a unit cube $(0, 1)^3$ and is uniformly partitioned into a Cartesian grid with $N$ intervals in each space direction.
A grid node is denoted as $(x_i,y_j,z_k)=(ih,jh,kh)$, $i,j,k = 1,2,\cdots,N$ where $h=1/N$ is the mesh size.

After time discretization with the DG or modified DG scheme, the heat equation can, after some rearrangement, be transformed into a sequence of 1D ODEs:
\begin{align}
    &\partial_{xx} u_{j,k} - \kappa^2 u_{j,k} = f_{j,k}, \quad \text{in } \Omega\cap\{y = y_j, z = z_k\}, \text{ for } j,k = 1, 2,\cdots, N-1,\\
    &\partial_{yy} u_{k,i} - \kappa^2 u_{k,i} = f_{k,i}, \quad \text{in } \Omega\cap\{z = z_k, x = x_i\}, \text{ for } k,i = 1, 2,\cdots, N-1,\\
    &\partial_{zz} u_{i,j} - \kappa^2 u_{i,j} = f_{i,j}, \quad \text{in } \Omega\cap\{x = x_i, y = y_j\}, \text{ for } i,j = 1, 2,\cdots, N-1.
\end{align}
subject to Dirichlet boundary conditions, respectively,
\begin{align}
    &u_{j,k}(\bm{x}) = g_D(\bm{x}, t^{n+1}), \quad \text{on } \partial\Omega\cap\{y = y_j, z = z_k\}, \text{ for } j,k = 1, 2,\cdots, N-1,\\
    &u_{k,i}(\bm{x}) = g_D(\bm{x}, t^{n+1}), \quad \text{on } \partial\Omega\cap\{z = z_k, x = x_i\}, \text{ for } k,i = 1, 2,\cdots, N-1,\\
    &u_{i,j}(\bm{x}) = g_D(\bm{x}, t^{n+1}), \quad \text{on } \partial\Omega\cap\{x = x_i, y = y_j\}, \text{ for } i,j = 1, 2,\cdots, N-1,
\end{align}
where $\kappa^2 = 2/\tau$, $u_{j,k}$, $u_{k,i}$, and $u_{i,j}$ are solutions to the 1D subproblems, and $f_{j,k}$, $f_{k,i}$, and $f_{i,j}$ are right-hand sides consisting of previous solutions.
Since the boundary correction technique cannot be applied here, the boundary conditions for the intermediate variables are taken directly from the physical boundary conditions at $t^{n+1}$.

In the 1D sub-problems, the domain $\Omega\cap\{y = y_j, z = z_k\}$ may consist of several disjoint subsets.
However, to illustrate the idea behind the 1D KFBI method, it suffices to consider the following simple two-point boundary value problem:
\begin{equation}\label{eqn:tp-bvp}
   \partial_{xx} u(x) - \kappa^2 u(x) = f(x), \quad \text{for }x\in(a, b)\subset (0, 1), 
\end{equation}
subject to the boundary condition
\begin{equation}
  u(x) = g(x), \quad x\in\{a, b\}.  
\end{equation}
Solving the 1D ODE is rather simple if the interval $(a,b)$ is partitioned into a fitted grid and standard finite difference methods or finite element methods are applied.
That is not the case here, since the endpoints of the one-dimensional intervals generally do not coincide with Cartesian grid nodes, leading to an unfitted mesh on the interval $(a,b)$.
In the KFBI method, the larger interval $(0, 1)$, which is uniformly partitioned into a grid, is the actual computational domain instead of $(a,b)$.
For each $x\in(0, 1)$, let $G(y,x)$ be the Green function such that it satisfies 
\begin{align}\label{eqn:green-1d}
    &\partial_{yy}G(y, x) - \kappa^2 G(y, x) = \delta(y-x), \quad y \in (0, 1),\\
    &G(0,x) = G(1,x) = 0,
\end{align}
where $\delta(\cdot)$ is the Dirac delta function in 1D.
As in higher-dimensional cases, the solution $u(x)$ of the boundary value problem \eqref{eqn:tp-bvp} can be expressed as the sum of a layer potential and a volume potential:
\begin{equation}
    u(x) = \partial_{y}G(y,x)\varphi(y)\big|_{y=a}^b + \int_a^b G(y, x)f(y) \,dy, \quad \text{for } x\in(0, 1),
\end{equation}
where $\partial_{y}G(y,x)\varphi(y)\big|_{y=a}^b =\partial_{y}G(b,x)\varphi(b) - \partial_{y}G(a,x)\varphi(a) $ is actually a 1D version of a boundary integral.
By incorporating the boundary condition, the above boundary value problem can be reduced to the boundary integral equation
\begin{equation}\label{eqn:bie-1d}
    \dfrac{1}{2}\varphi(x) + \partial_{y}G(y,x)\varphi(y)\big|_{y=a}^b = g(a) - \int_a^b G(y, x)f(y) \,dy,\quad \text{for }x\in\{a,b\}.
\end{equation}
The matrix-free GMRES method \cite{Saad2003}, which only requires the matrix-vector multiplication operator, is used to iteratively solve the discrete system of \eqref{eqn:bie-1d}.
In the KFBI method, boundary and volume integrals in \eqref{eqn:bie-1d} are all evaluated by solving equivalent but simple interface problems so that analytical expressions of integral kernels, consisting of the Green function and its derivatives, are never used.

Denote by $v(x) = \int_a^b G(y, x)f(y) \,dy$ the volume integral term. It satisfies the equivalent interface problem
\begin{equation}\label{eqn:ifp-v}
    \left\{
    \begin{aligned}
    &\partial_{xx}v - \kappa^2 v = \widetilde{f}(x), &\quad \text{ in }(0,1) \backslash \{a,b\},\\
    &[v](x) = 0, &\quad \text{ on }\{a,b\},\\
    &[\partial_x v](x) = 0, &\quad \text{ on }\{a,b\},\\
    &v(x) = 0, &\quad \text{ on }\{0,1\},
    \end{aligned}
    \right.
\end{equation}
where $\widetilde{f}(x)$ is an arbitrary extension of $f(x)$ to the larger interval $(0,1)$.
Here, the bracket $[\cdot]$ denotes the jump value of the one-sided limit from the interior of $(a,b)$ to the exterior.
The GMRES method starts with an initial guess $\varphi^0$ and generates a sequence of $\varphi^k, k = 1, 2, \cdots$ until it converges.
Let $w^k(x) = \partial_{y}G(y,x)\varphi^k(y)|_{y=a}^b$ be the double layer potential, which satisfies the equivalent interface problem
\begin{equation}\label{eqn:ifp-w}
    \left\{
    \begin{aligned}
    &\partial_{xx}w^k - \kappa^2 w^k = 0, &\quad \text{ in }(0,1) \backslash \{a,b\},\\
    &[w^k](x) = \varphi^k(x), &\quad \text{ on }\{a,b\},\\
    &[\partial_x w^k](x) = 0, &\quad \text{ on }\{a,b\},\\
    &w^k(x) = 0, &\quad \text{ on }\{0,1\}.
    \end{aligned}
    \right.
\end{equation}
From the interface problems, one can deduce jump values of the solution $w^k$ and $v$ and their derivatives as follows
\begin{align}
    &[v] = 0,\quad [\partial_x v] = 0, \quad [\partial_{xx} v] = [\widetilde{f}], \label{eqn:jumps-v}\\
    &[w^k] = \varphi^k,\quad [\partial_x w^k] = 0, \quad [\partial_{xx}w^k] = 0. \label{eqn:jumps-w}
\end{align}
The two interface problems are discretized with a second-order three-point finite difference method on the Cartesian grid over the interval $(0, 1)$.
In the vicinity of the interface, due to the low regularity of the solution, the local truncation error is large and leads to a poor approximation of the interface problem.
This difficulty can be mitigated by adding correction terms, computed from the jump conditions, to the right-hand side of the finite difference equations as defect corrections \cite{Zhou2023}.
As a result, the local truncation error becomes sufficiently small to obtain accurate approximations.
For example, let $v_i$ for $i=0,1,\cdots, N$ be the numerical approximation to the solution of \eqref{eqn:ifp-v} at the grid node $x_i$.
The corrected finite difference scheme is given by
\begin{equation}\label{eqn:crc-fds}
    \dfrac{v_{i+1}-2v_{i}+v_{i-1}}{h^2} - \kappa^2 v_i = \widetilde{f}(x_i) + C_i, \quad \text{for }i=1,2,\cdots, N-1,
\end{equation}
together with the boundary condition $v_0=v_N = 0$.
Here, $C_i$ is the correction term, which can be computed with the jump values in \eqref{eqn:jumps-v}.
We emphasize that the linear systems for the 1D interface problems remain tridiagonal because corrections are made only to the right-hand sides.
These systems can be solved with a highly efficient LU decomposition method, namely, the Thomas algorithm.
In addition, the tridiagonal linear systems are all identical if the problem has a constant coefficient and is discretized with the same mesh size in each spatial direction.
At the implementation level, one needs only to perform the LU decomposition in advance and store the results for subsequent computations.
Solving subsequent 1D interface problems only requires forward and backward substitutions.

The integral operators in the boundary integral equation \eqref{eqn:bie-1d} can be interpreted as one-sided limits of the two non-smooth potential functions $v(x)$ and $w^k(x)$ at the interface.
To be specific, we have
\begin{align}
    \int_a^b G(y, x)f(y) \,dy = v(x+) = v(x-), \quad \text{for }x\in\{a, b\}, \\
    \partial_{y}G(y,x)\varphi^k(y)|_{y=a}^b = \dfrac{1}{2}(w^k(x+) + w^k(x-)), \quad \text{for }x\in\{a, b\}.
\end{align}
Given the numerical solutions $v_i$ and $w_i^k$ for $i=0,1,\cdots,N$, let  $v_h(x)$ and $w_h^k(x)$ be two grid functions such that $v(x_i)=v_i$ and $w_h^k(x_i) = w_i^k$ for all $i=0,1,\cdots,N$.
One can compute the integral operators by interpolating from the grid functions $v_h(x)$ and $w^k_h(x)$.
Due to the jump of the potential functions on the interface $\{a,b\}$, one should also take into account the jump values in the Lagrangian interpolation.
More details on the correction and interpolation methods are described in \cite{Zhou2023}, and we omit them here.

\begin{remark}
When $\kappa$ is a constant, it is not too difficult to obtain an analytical expression of $G(y,x)$ and use it to form the $2\times 2$ system in \eqref{eqn:bie-1d}, which is then solved with a direct solver.
Once the unknown density $\varphi$ is determined, one only needs to solve an equivalent interface problem with known jump conditions, and the computational cost is comparable to that of solving a tridiagonal system.
In this case, the KFBI method is similar to Mayo's method \cite{Mayo1984, Mayo1985, Mayo1992}.
The KFBI method is more general than traditional boundary integral methods in that it can also handle variable-coefficient problems \cite{Ying2014}. We do not pursue that extension here.
\end{remark}

\subsection{Level set-ADI method for the Stefan problem}
\subsubsection{Level set method}
For moving-interface problems, the interface $\Gamma$ is time-dependent and must be determined together with the PDE in the bulk domain.
Here, the level set method is used to implicitly capture the position of the moving interface.
Given a level set function $\phi(\bm{x}):\mathcal{B}\rightarrow \mathbb{R}$, the interface $\Gamma$ is defined as the zero contour of $\phi(\bm{x})$:
\begin{equation}
    \Gamma = \left\{\bm{x}: \phi(\bm{x}) = 0, \bm{x}\in\mathcal{B}  \right\}.
\end{equation}
For a time-dependent interface, the time-dependent level set function $\phi(\bm{x}, t):\mathcal{B}\times[0,T]\rightarrow \mathbb{R}$ satisfies the Hamilton-Jacobi equation \cite{Osher2003}
\begin{equation}\label{eqn:hj-eqn}
    \phi_t + V_n |\nabla \phi| = 0.
\end{equation}
where $V_n$ is the prescribed normal velocity field of $\Gamma$.
Assume that there exists $\varepsilon_0 > 0$ and $\varepsilon_V \geq \varepsilon_0$.
A continuous extension of the normal velocity $V_n$ from $\Gamma$ to the domain $\mathcal{B}$ can be obtained by rewriting the Gibbs-Thomson relation \eqref{eqn:gibbs} as
\begin{equation}\label{eqn:vn-for}
    V_n = -\dfrac{\varepsilon_C(\phi)}{\varepsilon_V(\phi)} \nabla\cdot(\dfrac{\nabla\phi}{|\nabla\phi|}) -\dfrac{u}{\varepsilon_V(\phi)},
\end{equation}
Plugging \eqref{eqn:vn-for} into \eqref{eqn:hj-eqn} gives the governing equation for the evolution of $\phi$,
\begin{equation}\label{eqn:ls-for}
    \phi_t -\dfrac{u}{\varepsilon_V(\phi)} |\nabla \phi| = \dfrac{\varepsilon_C(\phi)}{\varepsilon_V(\phi)} \nabla\cdot(\dfrac{\nabla\phi}{|\nabla\phi|}) |\nabla \phi|.
\end{equation}
Note that the right-hand side of \eqref{eqn:ls-for} is a second-order term.
An explicit time stepping scheme requires a stability constraint on the time step $\tau<C h^2$ where $C$ is a constant depending on the ratio $\varepsilon_C/\varepsilon_V$.
The second-order stability constraint is induced by the mean curvature term and is stricter than the common Courant-Friedrichs-Lewy (CFL) condition that is encountered in level set methods.
To mitigate this problem, based on the relation $\nabla\cdot(\nabla\phi/|\nabla\phi|) |\nabla \phi| = \Delta\phi +  \nabla\cdot(1/|\nabla\phi|) |\nabla \phi|$, we use a semi-implicit level set method \cite{Smereka2003, Salac2008} by adding an extra term $S(\Delta\phi^{n+1} - \Delta \phi^n)$ to stabilize the explicit time-stepping scheme:
\begin{equation}\label{eqn:dis-ls}
    \dfrac{\phi^{n+1}-\phi^n}{\tau} -\dfrac{u^n}{\varepsilon_V(\phi^n)} |\nabla \phi^n| = \dfrac{\varepsilon_C(\phi^n)}{\varepsilon_V(\phi^n)} \nabla\cdot(\dfrac{\nabla\phi^n}{|\nabla\phi^n|+\epsilon}) |\nabla \phi^n| + S(\Delta\phi^{n+1} - \Delta \phi^n),
\end{equation}
where $S \geq |\varepsilon_C(\phi^n)/\varepsilon_V(\phi^n)|$ is a selectable stabilization coefficient and $\epsilon>0$ is a small constant to avoid division by zero.
Here, all but the stabilization term are treated explicitly to decouple the computation of $\phi^{n+1}$ and $u^{n+1}$ and avoid solving nonlinear systems.
The fifth-order WENO scheme \cite{Osher1991, Osher2003} is used for $|\nabla\phi^n|$ on the left-hand side of \eqref{eqn:dis-ls} since it represents convection in the normal direction.
Spatial derivatives on the right-hand side of \eqref{eqn:dis-ls} are approximated with central differences since they are diffusive in nature.

In practice, the level set function should be chosen as a signed distance function, which satisfies $|\nabla\phi|=1$, to improve the accuracy of the level set method.
As the level set function evolves, it gradually becomes too flat or too steep in some areas and fails to be a signed distance function.
To recover the signed-distance property, the level set function $\phi$ is reinitialized by solving the re-initialization equation \cite{Osher2003}
\begin{equation}\label{eqn:rein-eqn}
    \left\{
    \begin{aligned}
    &\dfrac{\partial\phi}{\partial\tau} + S(\phi_0) (|\nabla \phi|-1) = 0, \\
    &\phi(\bm{x}, 0) = \phi_0(\bm{x}).
    \end{aligned}
    \right.
\end{equation}
where $\tau$ is a pseudo-time, $\phi_0$ is the level set function before re-initialization, and $S(\phi) = \frac{\phi}{\sqrt{\phi^2 + h^2}}$ is a regularized sign function.
The fifth-order WENO scheme and the third-order TVD Runge-Kutta scheme are used to solve \eqref{eqn:rein-eqn}.

To improve computational efficiency, we implement a narrow-band level set method with a cutoff procedure, which leads to a local semi-implicit level set method \cite{Salac2008}.
The equations \eqref{eqn:ls-for} and \eqref{eqn:rein-eqn} are solved only in the narrow band rather than in the entire domain.
At each time step, the narrow band is computed following a layer-by-layer approach.
Define the first layer of irregular grid nodes by
\begin{equation}
    \mathcal{L}_1 = \{(x_i,y_j,z_k)|\exists(l, m,n)\in\mathcal{I}, \phi_{i,j,k}\phi_{i+l, j+m, k+n} < 0\},
\end{equation}
where the index set $\mathcal{I}$ is given as
\begin{equation}
    \mathcal{I} = \{(1,0,0), (-1,0,0), (0,1,0), (0,-1,0), (0,0,1), (0,0,-1)\}.
\end{equation}
Subsequently, the $k$-th layer of irregular grid nodes for $k=2, 3, \ldots$ are defined by
\begin{equation}
    \mathcal{L}_k = \{(x_i,y_j,z_k)|\exists(l, m,n)\in\mathcal{I}, (x_{i+l}, y_{j+m}, z_{k+n}) \in \mathcal{L}_{k-1} \}.
\end{equation}
Let the first $K$ layers of irregular grid nodes be the computational narrow band.
The grid nodes in the narrow band are called activated grid nodes, where level set function values are updated when solving \eqref{eqn:ls-for} and \eqref{eqn:rein-eqn}.

After the level set function has been updated at the activated nodes for evolution or re-initialization, it becomes discontinuous across the interface between activated and non-activated grid nodes, which may cause numerical oscillations and eventually contaminate the solution.
A simple cutoff procedure resolves this issue by modifying the level set function after reinitialization so that continuity is restored.
The cutoff function $\phi^{new} = c(\phi)$ is defined as
\begin{equation}
    c(\phi) = 
    \left\{
    \begin{aligned}
        &\phi,& \text{ if }   |\phi| < \phi^c,\\
        &\dfrac{\phi}{|\phi|}\phi^c, &\text{ if } |\phi| \geq \phi^c.
    \end{aligned}
    \right.
\end{equation}
where $\phi^c\in(0, Kh]$ is a cutoff threshold.
In this work, we choose $K = 8$ and $\phi^c=5h$.
Since the level set function is of little interest at grid nodes away from the zero contour, its value can be simply set as a constant to suppress numerical oscillations and simplify the algorithm as well.

\subsubsection{ADI scheme for the heat equation with an interface}
The evolution of the temperature field satisfies the heat equation and the Stefan condition, leading to an interface problem,
\begin{equation}\label{eqn:heat-ifp}
    \left\{
    \begin{aligned}
    &u_t = \Delta u, &\quad \text{in }\Omega_s\cup\Omega_l, \\
    &[u] = 0, &\quad \text{on }\Gamma,\\
    &[\partial_{\bm{n}}u] = -\phi_t/|\nabla\phi|,&\quad \text{on }\Gamma,
    \end{aligned}
    \right.
\end{equation}
where we have substituted $V_n$ with $-\phi_t/|\nabla\phi|$ by using \eqref{eqn:hj-eqn}.
The zeroth-order jump condition $[u] = 0$ is implied by the Gibbs-Thomson relation \eqref{eqn:gibbs} since the temperature is continuous across $\Gamma$.

When dimension splitting is carried out with ADI schemes for this problem, not only the PDE but also the interface condition must be decomposed into 1D forms.
Note that the zeroth-order jump condition $[u]=0$ can be reduced directly to a 1D condition.
To decompose the first-order jump condition, we take tangential derivatives of the zeroth-order jump condition to generate two additional first-order jump conditions.
Let $\boldsymbol{\tau}_1$ and $\boldsymbol{\tau}_2$ be two unit tangential vectors with $\boldsymbol{\tau}_1 \cdot\boldsymbol{\tau}_2=0$.
Taking tangential derivatives of both sides of $[u]=0$ in the directions of $\boldsymbol{\tau}_1$ and $\boldsymbol{\tau}_2$, and combining them with the first-order jump condition in \eqref{eqn:heat-ifp}, gives
\begin{equation}\label{eqn:jmp-sys}
    \boldsymbol{\tau}_1 \cdot [\nabla u] = 0,\quad \boldsymbol{\tau}_2 \cdot [\nabla u] = 0, \quad \bm{n}\cdot [\nabla u] = -\phi_t/|\nabla\phi| .
\end{equation}
From the $3\times 3$ linear system in \eqref{eqn:jmp-sys}, the first-order 1D jump conditions $[u_x]$, $[u_y]$, and $[u_z]$ can be obtained easily:
\begin{equation}
    [\nabla u] = -\phi_t/|\nabla\phi| \bm{n} = -\phi_t\nabla \phi /|\nabla\phi|^2,
\end{equation}
where we have used the relation $\bm{n} = \nabla \phi / |\nabla \phi|$.
At the discrete level, jump conditions are approximated by
\begin{equation}\label{eqn:jmp-dis}
    [u^{n+1}] = 0,\quad 
    [\nabla u^{n+1}] = \dfrac{\phi^{n} - \phi^{n+1}}{\tau }\dfrac{\nabla \phi^{n+1} }{|\nabla \phi^{n+1}|^2+\epsilon} := Q(\phi^{n}, \phi^{n+1}).
\end{equation}
where $\epsilon>0$ is also a small constant to avoid division by zero.

Based on the ADI schemes \eqref{eqn:DG-adi} and \eqref{eqn:3L-ADI}, the ADI scheme for the heat equation with an interface requires solving a sequence of 1D interface problems,
\begin{align}
    &\partial_{xx} u_{j,k} - \kappa^2 u_{j,k} = f_{j,k}, \quad \text{in } (\Omega_s\cup\Omega_l)\cap\{y = y_j, z = z_k\}, \text{ for } j,k = 1, 2,\cdots, N-1,\\
    &\partial_{yy} u_{k,i} - \kappa^2 u_{k,i} = f_{k,i}, \quad \text{in } (\Omega_s\cup\Omega_l)\cap\{z = z_k, x = x_i\}, \text{ for } k,i = 1, 2,\cdots, N-1,\\
    &\partial_{zz} u_{i,j} - \kappa^2 u_{i,j} = f_{i,j}, \quad \text{in } (\Omega_s\cup\Omega_l)\cap\{x = x_i, y = y_j\}, \text{ for } i,j = 1, 2,\cdots, N-1.
\end{align}
subject to interface conditions, respectively,
\begin{align}
    &[u_{j,k}] = 0, [\partial_x u_{j,k}] = Q(\phi^n,\phi^{n+1}) , \quad \text{on } \Gamma\cap\{y = y_j, z = z_k\}, \text{ for } j,k = 1, 2,\cdots, N-1, \\
    &[u_{k,i}] = 0, [\partial_y u_{k,i}] = Q(\phi^n,\phi^{n+1}) , \quad \text{on } \Gamma\cap\{z = z_k, x = x_i\}, \text{ for } k,i = 1, 2,\cdots, N-1,\\
    &[u_{i,j}] = 0, [\partial_z u_{i,j}] = Q(\phi^n,\phi^{n+1}) , \quad \text{on } \Gamma\cap\{x = x_i, y = y_j\}, \text{ for } i,j = 1, 2,\cdots, N-1.
\end{align}

The 1D interface problems are similar to those in \eqref{eqn:ifp-v} and \eqref{eqn:ifp-w} and can be solved with the finite difference scheme in \eqref{eqn:crc-fds} and the fast Thomas algorithm.
Unlike the case of solving boundary value problems, there is no need to solve boundary integral equations since jump conditions of the 1D interface problems are already known.
The computational cost of solving the interface problems is almost the same as in cases without interfaces.

\section{Numerical results}\label{sec:result}

\begin{figure}[htbp]
    \centering
    \subfigure[]{\includegraphics[width=0.4\linewidth]{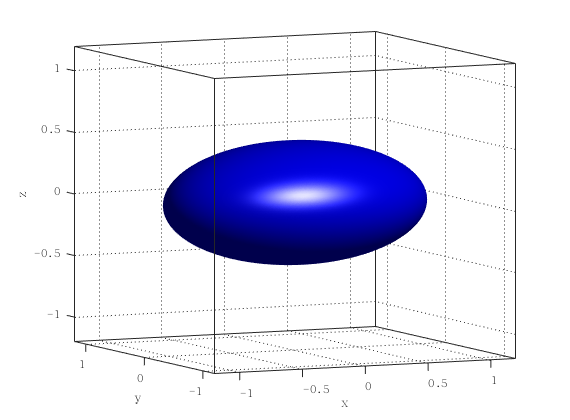}}
    \subfigure[]{\includegraphics[width=0.4\linewidth]{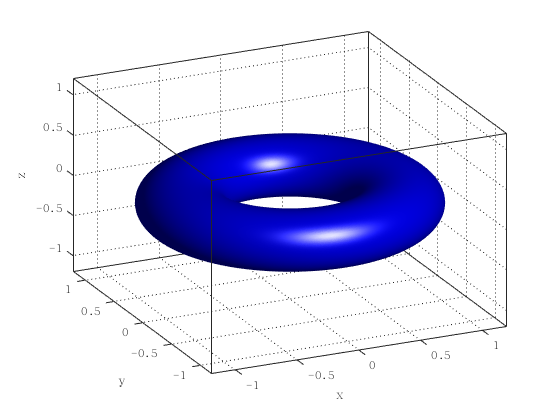}}\\
    \subfigure[]{\includegraphics[width=0.4\linewidth]{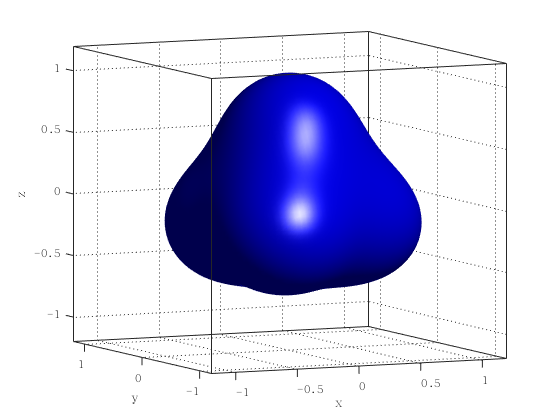}}
    \subfigure[]{\includegraphics[width=0.4\linewidth]{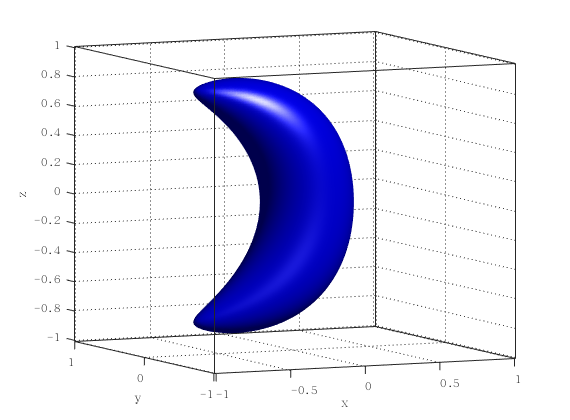}}
    \caption{Irregular domains: (a) an ellipsoid; (b) a torus; (c) a four-atom molecule; (d) a banana.}
    \label{fig:domains}
\end{figure}

In this section, we present numerical examples for the heat equation and a reaction-diffusion equation on fixed domains, as well as for the Stefan problem with a free boundary, to assess the performance of the proposed method.
The following level set functions $\phi(\bm{x})$ for defining irregular domains are used in the numerical tests:
\begin{itemize}
    \item An ellipsoid: $$\phi(\bm{x}) = \frac{x^{2}}{1^{2}} + \frac{y^{2}}{0.7^{2}} + \frac{z^{2}}{0.5^{2}} - 1; $$
    \item A torus: $$ \phi(\bm{x}) = (\sqrt{x^{2} + y^{2}} - 0.8)^{2} + z^{2} - 0.34^2; $$
    \item A four-atom molecule: $$\phi(\bm{x}) = c - \sum_{i = 1}^{4} \exp(-\frac{||\bm{x}-\bm{x}_i||^2}{r^{2}}),$$
    with $c=r=0.6$ and $\bm{x}_1 = (\frac{\sqrt{3}}{3}, 0, -\frac{\sqrt{6}}{12})$, $\bm{x}_2 = (-\frac{\sqrt{3}}{6}, 0.5, - \frac{\sqrt{6}}{12})$, $\bm{x}_3 = (-\frac{\sqrt{3}}{6}, -0.5, -\frac{\sqrt{6}}{12})$, $\bm{x}_4 = (0, 0, \frac{\sqrt{6}}{4})$;
    \item A banana: 
    \begin{align*}
         \phi(\bm{x}) &= (7x + 6)^{4} + 2401y^{4} + 3601.5z^{4} + 98(7x+6)^{2}y^{2} +98(7x+6)^{2}z^{2} \\ & +4802y^{2}z^{2}-94(7x+6)^{2}+3822y^{2}-4606z^{2}+1521.
    \end{align*}
\end{itemize}
The domain boundaries $\Gamma=\partial\Omega$ defined by the zero contours of the level set functions $\phi(\bm{x})$ are plotted in \Cref{fig:domains}.

Numerical errors are estimated on the grid nodes in $\Omega^h$ at the final time $t=T$ in both the $L_2$ and $L_{\infty}$ norms, which are defined as
\begin{equation}
    \Vert e^h\Vert_{L_2} = \sqrt{\dfrac{1}{N_{\Omega^h}}\sum_{\boldsymbol{x}\in\Omega^h}|u^h(\boldsymbol{x}, T) - u(\boldsymbol{x}, T)|^2}, \quad \Vert e^h\Vert_{L_{\infty}} = \max_{\boldsymbol{x}\in\Omega^h}|u^h(\boldsymbol{x}, T)-u(\boldsymbol{x},T)|,
\end{equation}
where $N_{\Omega^h}$ is the number of grid nodes in $\Omega^h$, and $u^h$ and $u$ are the numerical and exact solutions, respectively.
The convergence order is computed by
\begin{equation}
    \text{order} = \dfrac{\log (||e^{2h}||/||e^h||)}{\log 2}.
\end{equation}
All numerical experiments are performed on a personal computer with a 16-core Intel(R) Core(TM) i7-10700K CPU @ 3.80GHz. 
The codes for the numerical experiments are written in C++.

\subsection{The heat equation}
In the first example, an initial-boundary value problem (IBVP) for the heat equation \eqref{eqn:heat1} is solved:
\begin{equation}\label{eqn:heat1}
    u_{t}(\bm{x}, t) = a\Delta u(\bm{x}, t) + f(\bm{x}, t), 
\end{equation}
where the diffusion coefficient is chosen as $a = 2$.
Here, the initial and boundary conditions for $u(\bm{x}, t)$ and the source term $f(\bm{x}, t)$ are chosen such that the problem has the exact solution
\begin{equation}
    u(\bm{x}, t) = \sin\paren{\frac{1}{\sqrt{3}}(x + y + z) - t}.
\end{equation}
Two different shapes of $\Omega$, an ellipsoid and a torus, are considered in this example.
The bounding box is chosen as $\mathcal{B} = (-1.2,1.2)^3$.
The final time is set as $T=0.5$.

\subsubsection{Convergence test}
\begin{figure}[htbp]
    \centering
    \includegraphics[width=0.45\linewidth]{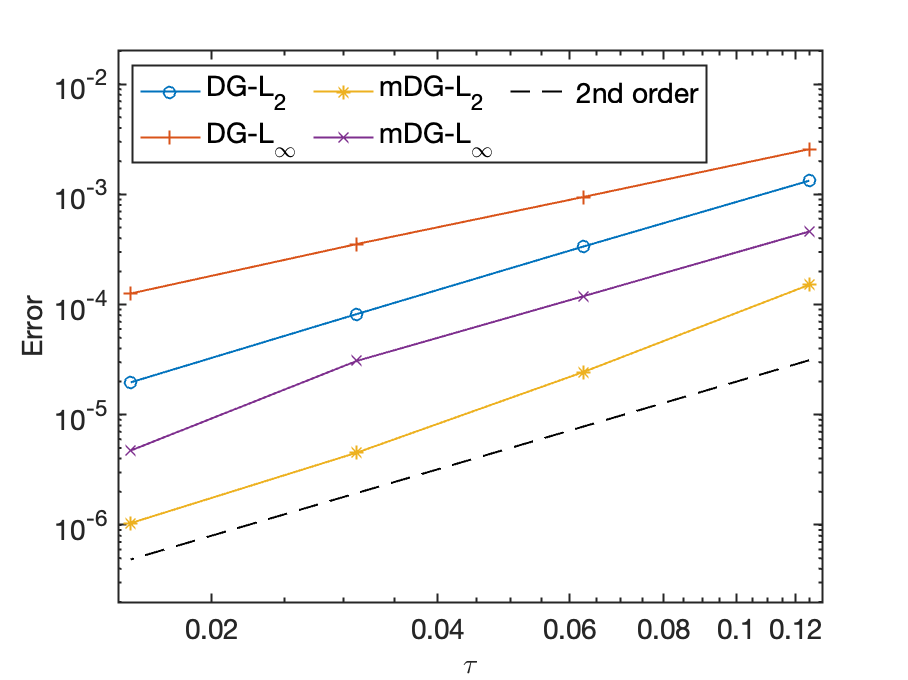}
    \includegraphics[width=0.45\linewidth]{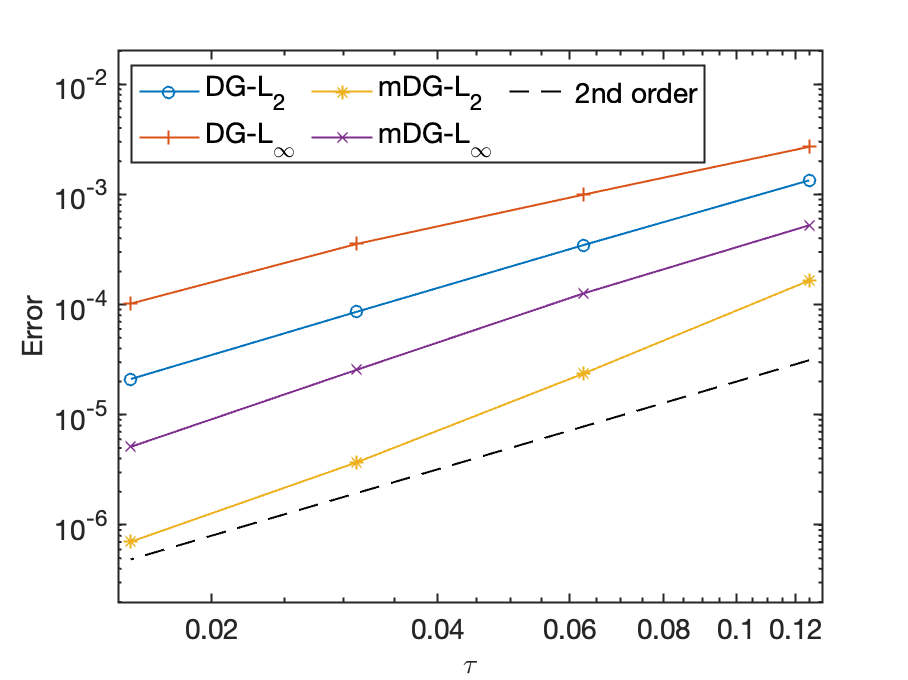}
    \caption{Temporal convergence of the two KFBI-ADI schemes for the heat equation on the ellipsoid-shaped domain (left) and the torus-shaped domain (right).}
    \label{fig:heat-dt}
\end{figure}

In order to study the overall convergence rates of the proposed methods, the spatial grid and time step are simultaneously refined with $h=2.4/N$ and $\tau=2/N$ to perform computations.
The numerical errors and convergence orders of the KFBI-ADI method based on the DG scheme and the modified DG scheme are summarized in \Cref{tab:ell-heat,tab:torus-heat}.
\begin{table}[htbp]
\centering
\caption{Numerical results of the heat equation on an ellipsoid-shaped domain.}
\label{tab:ell-heat}
\begin{tabular}{|c|cccc|cccc|}
\hline
\multirow{3}{*}{N} &
  \multicolumn{4}{c|}{KFBI-DG-ADI} &
  \multicolumn{4}{c|}{KFBI-mDG-ADI} \\ \cline{2-9} 
 &
  \multicolumn{2}{c|}{$L_2$} &
  \multicolumn{2}{c|}{$L_{\infty}$} &
  \multicolumn{2}{c|}{$L_2$} &
  \multicolumn{2}{c|}{$L_{\infty}$} \\ \cline{2-9} 
 &
  \multicolumn{1}{c|}{Error} &
  \multicolumn{1}{c|}{Rate} &
  \multicolumn{1}{c|}{Error} &
  Rate &
  \multicolumn{1}{c|}{Error} &
  \multicolumn{1}{c|}{Rate} &
  \multicolumn{1}{c|}{Error} &
  Rate \\ \hline
16 &
  \multicolumn{1}{c|}{1.32e-3} &
  \multicolumn{1}{c|}{} &
  \multicolumn{1}{c|}{2.35e-3} &
   &
  \multicolumn{1}{c|}{1.62e-4} &
  \multicolumn{1}{c|}{} &
  \multicolumn{1}{c|}{4.37e-4} &
   \\ \hline
32 &
  \multicolumn{1}{c|}{3.32e-4} &
  \multicolumn{1}{c|}{1.99} &
  \multicolumn{1}{c|}{7.74e-4} &
  1.60 &
  \multicolumn{1}{c|}{2.64e-5} &
  \multicolumn{1}{c|}{2.62} &
  \multicolumn{1}{c|}{1.15e-4} &
  1.93 \\ \hline
64 &
  \multicolumn{1}{c|}{8.14e-5} &
  \multicolumn{1}{c|}{2.03} &
  \multicolumn{1}{c|}{3.39e-4} &
  1.19 &
  \multicolumn{1}{c|}{4.93e-6} &
  \multicolumn{1}{c|}{2.42} &
  \multicolumn{1}{c|}{2.26e-5} &
  2.35 \\ \hline
128 &
  \multicolumn{1}{c|}{1.97e-5} &
  \multicolumn{1}{c|}{2.05} &
  \multicolumn{1}{c|}{1.26e-4} &
  1.43 &
  \multicolumn{1}{c|}{1.04e-6} &
  \multicolumn{1}{c|}{2.25} &
  \multicolumn{1}{c|}{4.74e-6} &
  2.25 \\ \hline
256 &
  \multicolumn{1}{c|}{4.77e-6} &
  \multicolumn{1}{c|}{2.05} &
  \multicolumn{1}{c|}{3.39e-5} &
  1.89 &
  \multicolumn{1}{c|}{2.39e-7} &
  \multicolumn{1}{c|}{2.12} &
  \multicolumn{1}{c|}{9.74e-7} &
  2.28 \\ \hline
\end{tabular}
\end{table}
\begin{table}[htbp]
\centering
\caption{Numerical results of the heat equation on a torus-shaped domain.}
\label{tab:torus-heat}
\begin{tabular}{|c|cccc|cccc|}
\hline
\multirow{3}{*}{N} &
  \multicolumn{4}{c|}{KFBI-DG-ADI} &
  \multicolumn{4}{c|}{KFBI-mDG-ADI} \\ \cline{2-9} 
 &
  \multicolumn{2}{c|}{$L_2$} &
  \multicolumn{2}{c|}{$L_{\infty}$} &
  \multicolumn{2}{c|}{$L_2$} &
  \multicolumn{2}{c|}{$L_{\infty}$} \\ \cline{2-9} 
 &
  \multicolumn{1}{c|}{Error} &
  \multicolumn{1}{c|}{Rate} &
  \multicolumn{1}{c|}{Error} &
  Rate &
  \multicolumn{1}{c|}{Error} &
  \multicolumn{1}{c|}{Rate} &
  \multicolumn{1}{c|}{Error} &
  Rate \\ \hline
16 &
  \multicolumn{1}{c|}{1.32e-3} &
  \multicolumn{1}{c|}{} &
  \multicolumn{1}{c|}{2.59e-3} &
   &
  \multicolumn{1}{c|}{1.75e-4} &
  \multicolumn{1}{c|}{} &
  \multicolumn{1}{c|}{5.12e-4} &
   \\ \hline
32 &
  \multicolumn{1}{c|}{3.48e-4} &
  \multicolumn{1}{c|}{1.92} &
  \multicolumn{1}{c|}{8.35e-4} &
  1.63 &
  \multicolumn{1}{c|}{2.49e-5} &
  \multicolumn{1}{c|}{2.81} &
  \multicolumn{1}{c|}{1.21e-4} &
  2.08 \\ \hline
64 &
  \multicolumn{1}{c|}{8.58e-5} &
  \multicolumn{1}{c|}{2.02} &
  \multicolumn{1}{c|}{3.41e-4} &
  1.29 &
  \multicolumn{1}{c|}{3.96e-6} &
  \multicolumn{1}{c|}{2.65} &
  \multicolumn{1}{c|}{2.52e-5} &
  2.26 \\ \hline
128 &
  \multicolumn{1}{c|}{2.11e-5} &
  \multicolumn{1}{c|}{2.02} &
  \multicolumn{1}{c|}{1.02e-4} &
  1.74 &
  \multicolumn{1}{c|}{7.04e-7} &
  \multicolumn{1}{c|}{2.49} &
  \multicolumn{1}{c|}{5.12e-6} &
  2.30 \\ \hline
256 &
  \multicolumn{1}{c|}{5.15e-6} &
  \multicolumn{1}{c|}{2.03} &
  \multicolumn{1}{c|}{3.52e-5} &
  1.53 &
  \multicolumn{1}{c|}{1.43e-7} &
  \multicolumn{1}{c|}{2.30} &
  \multicolumn{1}{c|}{1.18e-6} &
  2.12 \\ \hline
\end{tabular}
\end{table}
The difference in temporal convergence rates of the two schemes is also studied using a fine spatial grid ($N=256$) and successively refining the time step.
The results are plotted in \Cref{fig:heat-dt}.
The modified DG scheme produces more accurate results than the original DG scheme on both the ellipsoid and torus domains.
Both ADI schemes have approximately second-order accuracy.
However, the DG scheme exhibits order reduction in the $L_{\infty}$ norm.
The modified DG scheme is free of this problem since it has higher consistency of boundary conditions by construction. 
In addition, for all time steps considered here, which are much larger than those used in explicit schemes, the computation remains stable.
These results provide numerical evidence that the proposed ADI schemes are also unconditionally stable on irregular domains.

\subsubsection{Efficiency test}
An important advantage of ADI schemes is their efficiency for time-dependent problems in multiple spatial dimensions.
Since the method is based on dimension splitting, it is naturally well suited to parallelization.
In each sub-step of the ADI schemes, a sequence of one-dimensional subproblems, which are independent of each other, can be solved separately using simple loops.
In C++, custom implementations of the ADI schemes written with ordinary loops can be converted easily into a multithreaded version using the OpenMP library \cite{chandra2001parallel} by adding only a few directives.
Even without parallelization techniques, the ADI scheme with the fast Thomas algorithm is still very efficient due to its linear complexity with a small constant factor.

\begin{figure}[htbp]
    \centering
    \includegraphics[width=0.45\textwidth]{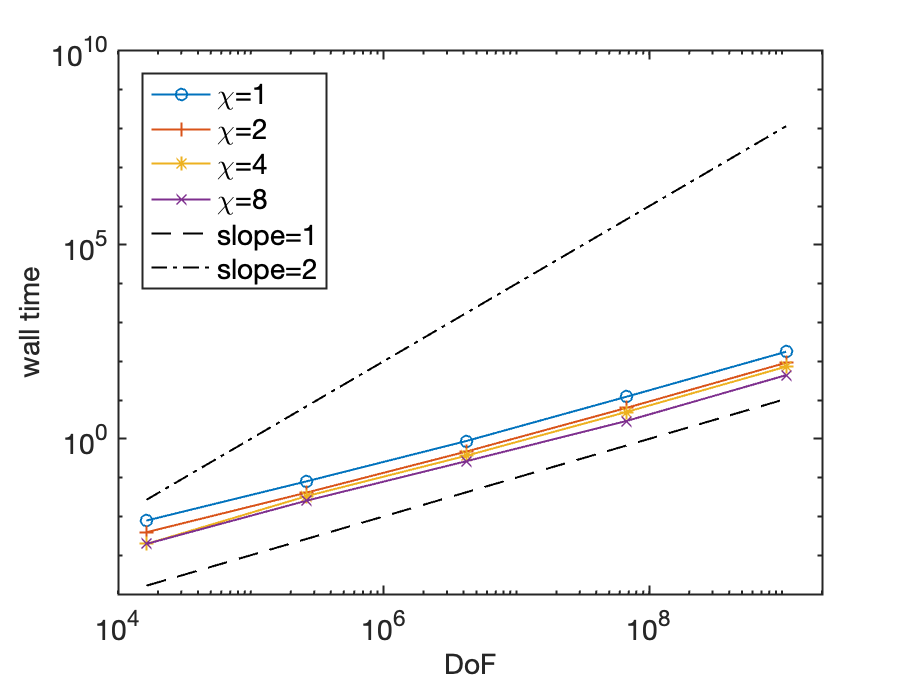}
    \includegraphics[width=0.45\textwidth]{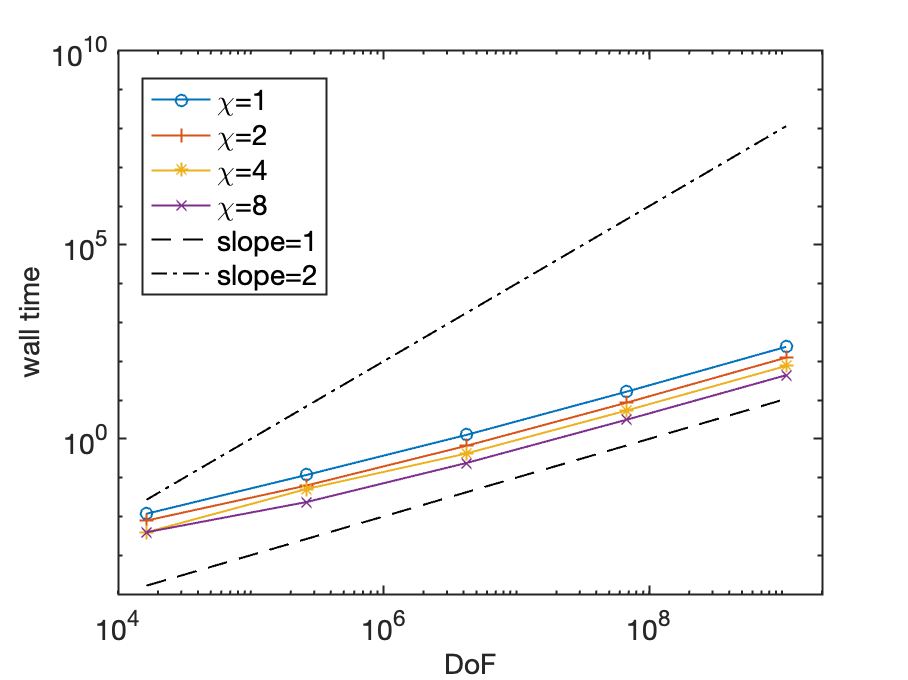}
    \caption{Wall times of the KFBI-ADI method with different mesh sizes and thread numbers for the heat equation on the ellipsoid-shaped domain (left) and the torus-shaped domain (right).}
    \label{fig:time1}
\end{figure}

To demonstrate the efficiency of the proposed methods and the speed-up obtained through multithreading, we solve the heat equation with the same configurations as before and collect wall times for different thread counts $\chi = 1, 2, 4,$ and $8$.
The results are plotted in \Cref{fig:time1}, in which the x-axis represents the number of degrees of freedom, defined as the total number of Cartesian grid nodes times the number of time steps, i.e. $DoF = (N+1)^3N_t$, and the y-axis represents the total wall time (in seconds) for numerical tests with different thread counts.
As the number of degrees of freedom increases, the slopes of these lines approach $1$ for different thread counts, indicating linear complexity of the proposed method.

Furthermore, detailed speed-up ratios are presented in \Cref{tab:speedup} to estimate the acceleration achieved by multithreading in the proposed method. 
Here, the speed-up ratio is defined as the ratio of wall times using one and $\chi$ threads, i.e.,
\begin{equation}
    \text{Speed-up ratio} = \dfrac{\text{wall time using one thread}}{\text{wall time using }\chi\text{ threads}}.
\end{equation}
From the results, it is clear that the method benefits substantially from multithreading, which shows that the proposed method is well suited for parallelization.
The speed-up ratio decreases as the number of threads increases and it fails to achieve the ideal case.
This is mostly due to the increase in communication and synchronization costs when more threads are used.

\begin{table}[htbp]
\centering
\caption{Speed-up ratios of the proposed method with different thread numbers.}
\label{tab:speedup}
\begin{tabular}{|c|ccc|ccc|}
\hline
    & \multicolumn{3}{c|}{ellipsoid domain}                               & \multicolumn{3}{c|}{torus domain}                                   \\ \hline
N & \multicolumn{1}{c|}{$\chi=2$} & \multicolumn{1}{c|}{$\chi=4$} & $\chi=8$ & \multicolumn{1}{c|}{$\chi=2$} & \multicolumn{1}{c|}{$\chi=4$} & $\chi=8$ \\ \hline
64  & \multicolumn{1}{c|}{1.86} & \multicolumn{1}{c|}{2.39} & 3.31 & \multicolumn{1}{c|}{1.91} & \multicolumn{1}{c|}{3.02} & 5.26 \\ \hline
128 & \multicolumn{1}{c|}{1.94} & \multicolumn{1}{c|}{2.47} & 4.22 & \multicolumn{1}{c|}{1.92} & \multicolumn{1}{c|}{3.06} & 5.29 \\ \hline
256 & \multicolumn{1}{c|}{1.92} & \multicolumn{1}{c|}{2.42} & 4.01 & \multicolumn{1}{c|}{1.90} & \multicolumn{1}{c|}{3.12} & 5.39 \\ \hline
\end{tabular}
\end{table}

\subsection{The reaction-diffusion equation}
In this example, an IBVP for a reaction-diffusion equation is solved with the proposed method.
We consider the Fisher equation
\begin{equation}\label{diffusion-reaction}
    u_{t} = \epsilon \Delta u + f(u), 
\end{equation}
where $\epsilon>0$ is a small parameter, assumed to be $0.1$, and $f(u) = \frac{2}{3 \epsilon}u^{2}(1 - u)$ is a reaction term.
The initial and boundary conditions are chosen so that the exact solution is
\begin{equation}
    u(\bm{x},t) = \dfrac{1}{1 + e^{(\bm{n}\cdot\bm{x}-t)/(3\epsilon)}},
\end{equation}
where $\bm{n}\in\mathbb{R}^3$ is a vector of unit length, representing the propagation direction of the reaction wavefront.
The problem is solved on a banana-shaped domain and a four-atom molecular domain.
The bounding box for the banana-shaped domain is chosen as $(-1.2,1.2)^3$, while that for the molecular domain is chosen as $(-1.01, 1.01)^3$.
The final time is set as $T=0.5$.
The scalar nonlinear algebraic equation resulting from the implicit step for the reaction term is solved by Newton's method with absolute tolerance $tol = 10^{-12}$.

Numerical errors and convergence orders obtained by two KFBI-ADI schemes are summarized in \Cref{tab:banana-rd,tab:molecular-rd}.
Similarly, both schemes exhibit second-order accuracy.
Since the exact solution in this example is nearly constant away from the moving front, the original DG scheme produces more accurate results than in the case with time-dependent boundary conditions and has accuracy comparable to that of the modified DG scheme.

\begin{figure}[htbp]
    \centering
    \subfigure[t=0]{\includegraphics[width=0.3\textwidth]{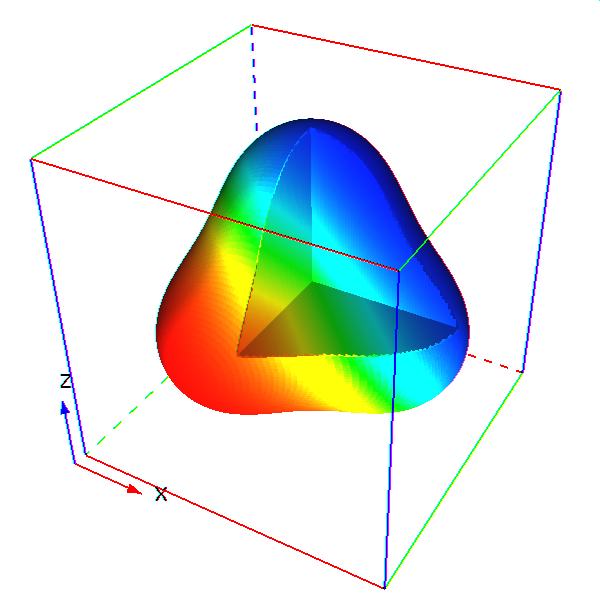}}
    \subfigure[t=0.25]{\includegraphics[width=0.3\textwidth]{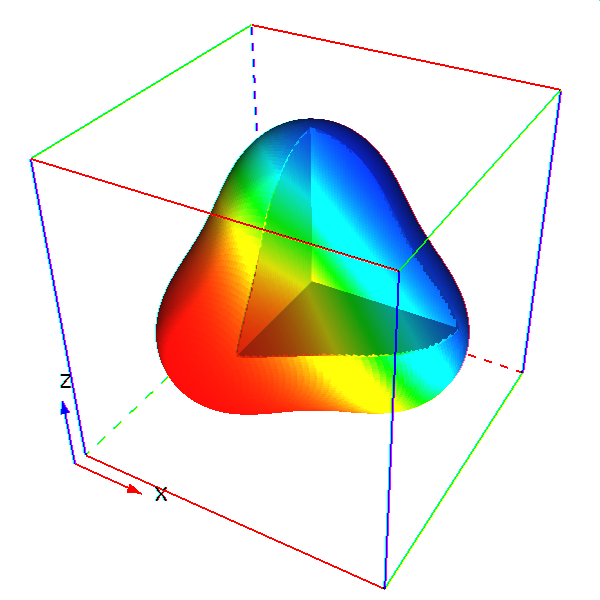}}
    \subfigure[t=0.5]{\includegraphics[width=0.3\textwidth]{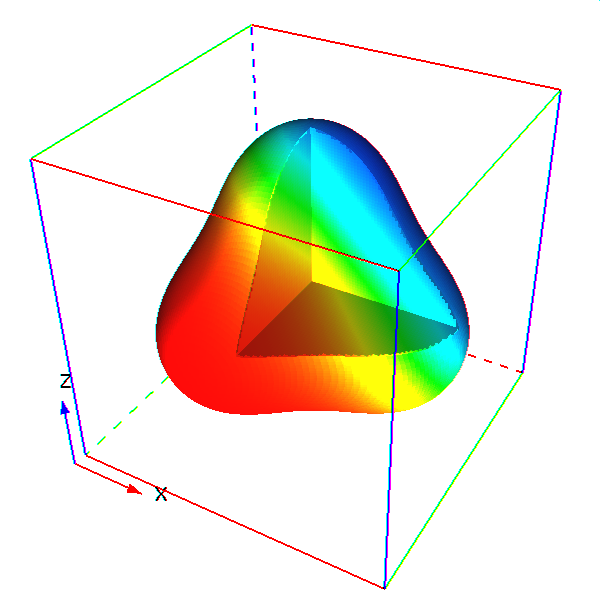}}
\caption{Numerical solution of the reaction-diffusion equation on a four-atom molecular domain.}
    \label{fig:molecular-rd}
\end{figure}

\begin{figure}[htbp]
    \centering
    \subfigure[t=0]{\includegraphics[width=0.3\textwidth]{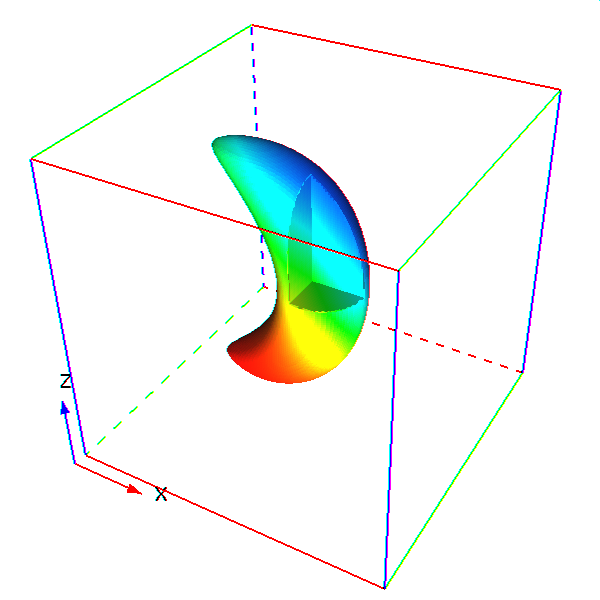}}
    \subfigure[t=0.25]{\includegraphics[width=0.3\textwidth]{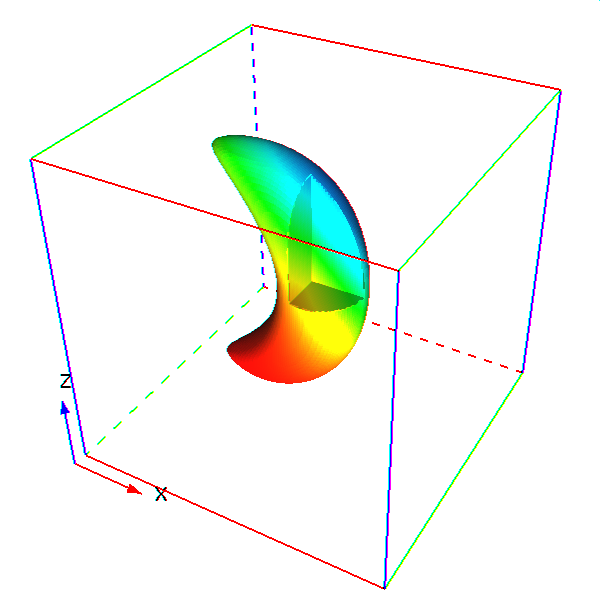}}
    \subfigure[t=0.5]{\includegraphics[width=0.3\textwidth]{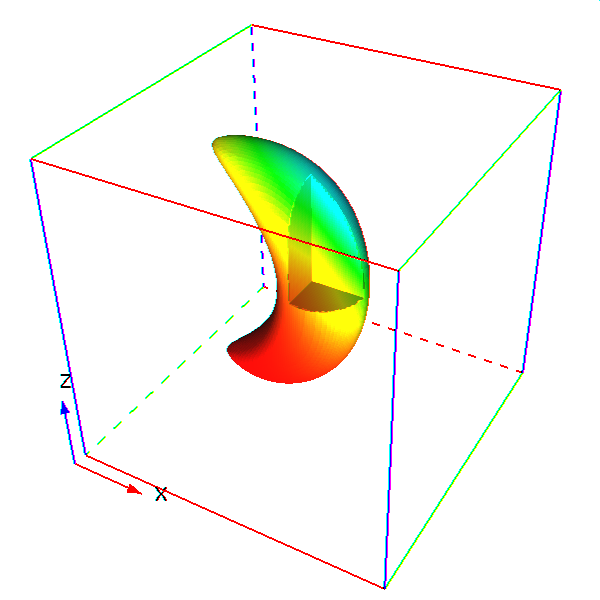}}
    \caption{Numerical solution to the reaction-diffusion equation on a banana-shaped domain.}
    \label{fig:ban-rd}
\end{figure}

The banana-shaped domain is singular and has two sharp corners, in contrast to the previous test domains.
The proposed ADI methods can naturally handle sharp corners without losing accuracy because the geometric singularity is no longer problematic in the 1D subproblems.
Snapshots of the numerical solutions at $t=0$, $t=0.25$, and $t=0.5$ are also displayed in \Cref{fig:molecular-rd,fig:ban-rd}, showing the time evolution of the moving front between the red and blue regions.

\begin{table}[htbp]
\centering
\caption{Numerical results of the reaction-diffusion equation on a four-atom molecular domain.}
\label{tab:molecular-rd}
\begin{tabular}{|c|cccc|cccc|}
\hline
\multirow{3}{*}{N} &
  \multicolumn{4}{c|}{KFBI-DG-ADI} &
  \multicolumn{4}{c|}{KFBI-mDG-ADI} \\ \cline{2-9} 
 &
  \multicolumn{2}{c|}{$L_2$} &
  \multicolumn{2}{c|}{$L_{\infty}$} &
  \multicolumn{2}{c|}{$L_2$} &
  \multicolumn{2}{c|}{$L_{\infty}$} \\ \cline{2-9} 
 &
  \multicolumn{1}{c|}{Error} &
  \multicolumn{1}{c|}{Rate} &
  \multicolumn{1}{c|}{Error} &
  Rate &
  \multicolumn{1}{c|}{Error} &
  \multicolumn{1}{c|}{Rate} &
  \multicolumn{1}{c|}{Error} &
  Rate \\ \hline
16 &
  \multicolumn{1}{c|}{1.24e-3} &
  \multicolumn{1}{c|}{} &
  \multicolumn{1}{c|}{5.63e-3} &
   &
  \multicolumn{1}{c|}{3.91e-3} &
  \multicolumn{1}{c|}{} &
  \multicolumn{1}{c|}{9.32e-3} &
   \\ \hline
32 &
  \multicolumn{1}{c|}{3.55e-4} &
  \multicolumn{1}{c|}{1.80} &
  \multicolumn{1}{c|}{1.89e-3} &
  1.57 &
  \multicolumn{1}{c|}{9.17e-4} &
  \multicolumn{1}{c|}{2.09} &
  \multicolumn{1}{c|}{2.32e-3} &
  2.01 \\ \hline
64 &
  \multicolumn{1}{c|}{9.34e-5} &
  \multicolumn{1}{c|}{1.93} &
  \multicolumn{1}{c|}{5.21e-4} &
  1.86 &
  \multicolumn{1}{c|}{2.20e-4} &
  \multicolumn{1}{c|}{2.06} &
  \multicolumn{1}{c|}{5.57e-4} &
  2.06 \\ \hline
128 &
  \multicolumn{1}{c|}{2.39e-5} &
  \multicolumn{1}{c|}{1.97} &
  \multicolumn{1}{c|}{1.38e-4} &
  1.92 &
  \multicolumn{1}{c|}{5.40e-5} &
  \multicolumn{1}{c|}{2.03} &
  \multicolumn{1}{c|}{1.36e-4} &
  2.03 \\ \hline
256 &
  \multicolumn{1}{c|}{6.03e-6} &
  \multicolumn{1}{c|}{1.99} &
  \multicolumn{1}{c|}{3.46e-5} &
  2.00 &
  \multicolumn{1}{c|}{1.34e-5} &
  \multicolumn{1}{c|}{2.01} &
  \multicolumn{1}{c|}{3.34e-5} &
  2.03 \\ \hline
\end{tabular}
\end{table}

\begin{table}[htbp]
\centering
\caption{Numerical results of the reaction-diffusion equation on a banana-shaped domain.}
\label{tab:banana-rd}
\begin{tabular}{|c|cccc|cccc|}
\hline
\multirow{3}{*}{N} &
  \multicolumn{4}{c|}{KFBI-DG-ADI} &
  \multicolumn{4}{c|}{KFBI-mDG-ADI} \\ \cline{2-9} 
 &
  \multicolumn{2}{c|}{$L_2$} &
  \multicolumn{2}{c|}{$L_{\infty}$} &
  \multicolumn{2}{c|}{$L_2$} &
  \multicolumn{2}{c|}{$L_{\infty}$} \\ \cline{2-9} 
 &
  \multicolumn{1}{c|}{Error} &
  \multicolumn{1}{c|}{Rate} &
  \multicolumn{1}{c|}{Error} &
  Rate &
  \multicolumn{1}{c|}{Error} &
  \multicolumn{1}{c|}{Rate} &
  \multicolumn{1}{c|}{Error} &
  Rate \\ \hline
16 &
  \multicolumn{1}{c|}{1.51e-3} &
  \multicolumn{1}{c|}{} &
  \multicolumn{1}{c|}{5.04e-3} &
   &
  \multicolumn{1}{c|}{2.42e-3} &
  \multicolumn{1}{c|}{} &
  \multicolumn{1}{c|}{6.63e-3} &
   \\ \hline
32 &
  \multicolumn{1}{c|}{4.74e-4} &
  \multicolumn{1}{c|}{1.67} &
  \multicolumn{1}{c|}{1.64e-3} &
  1.62 &
  \multicolumn{1}{c|}{5.58e-4} &
  \multicolumn{1}{c|}{2.12} &
  \multicolumn{1}{c|}{1.63e-3} &
  2.02 \\ \hline
64 &
  \multicolumn{1}{c|}{1.27e-4} &
  \multicolumn{1}{c|}{1.90} &
  \multicolumn{1}{c|}{4.78e-4} &
  1.78 &
  \multicolumn{1}{c|}{1.35e-4} &
  \multicolumn{1}{c|}{2.05} &
  \multicolumn{1}{c|}{4.10e-4} &
  1.99 \\ \hline
128 &
  \multicolumn{1}{c|}{3.28e-5} &
  \multicolumn{1}{c|}{1.95} &
  \multicolumn{1}{c|}{1.24e-4} &
  1.95 &
  \multicolumn{1}{c|}{3.31e-5} &
  \multicolumn{1}{c|}{2.03} &
  \multicolumn{1}{c|}{1.02e-4} &
  2.01 \\ \hline
256 &
  \multicolumn{1}{c|}{8.31e-6} &
  \multicolumn{1}{c|}{1.98} &
  \multicolumn{1}{c|}{3.21e-5} &
  1.95 &
  \multicolumn{1}{c|}{8.22e-6} &
  \multicolumn{1}{c|}{2.01} &
  \multicolumn{1}{c|}{2.53e-5} &
  2.01 \\ \hline
\end{tabular}
\end{table}

\subsection{The Stefan problem}
In the final example, we solve the Stefan problem to simulate dendritic solidification.
In this case, the free boundary is not only irregular but also time-dependent.
The level set method is used to capture the location of the free boundary.
A uniform Cartesian grid with $128\times 128\times 128$ cells over the bounding box $\mathcal{B} = (-2,2)^{3}$ is used for computations. 
The time step is determined from a CFL number of $0.5$ for the level set equation.
The computation is terminated when the free boundary reaches the box boundary.

Initially, a spherical solid seed with zero temperature is placed at the origin and is surrounded by undercooled liquid with a lower temperature; the Stefan number is $St = -0.5$.
Accordingly, the initial values of the level set function and the temperature field are chosen as
\begin{equation}
    \phi(\bm{x},0)=\phi_0(\bm{x}) = ||\bm{x}||-r, \quad T(\bm{x},0) = T_0(\bm{x}) = H_{\epsilon}(\phi(\bm{x})) St,
\end{equation}
where $r=0.1$ is the seed radius and $H_{\epsilon}(r) = \frac{1}{2}(1 + \tanh(r/\epsilon))$ is a regularized Heaviside function; here, $h$ is the mesh size. 
Anisotropic coefficients are chosen as 
\begin{equation}
    \varepsilon_C(\bm{n}) = \varepsilon_V(\bm{n}) = \bar{\varepsilon}(1-A(n_1^4+n_2^4+n_3^4)),
\end{equation}
where $n_1$, $n_2$ and $n_3$ are components of the unit normal vector $\bm{n}$.
Surface tension and molecular kinetic effects stabilize the free boundary by suppressing the unstable growth of small perturbations. 
Due to the anisotropy, the spherical solid has preferred growth directions corresponding to vectors $\bm{n}$ for which $\varepsilon_C(\bm{n})$ and $\varepsilon_V(\bm{n})$ are small.

\begin{figure}[htbp]
    \centering
    \includegraphics[width=0.8\textwidth]{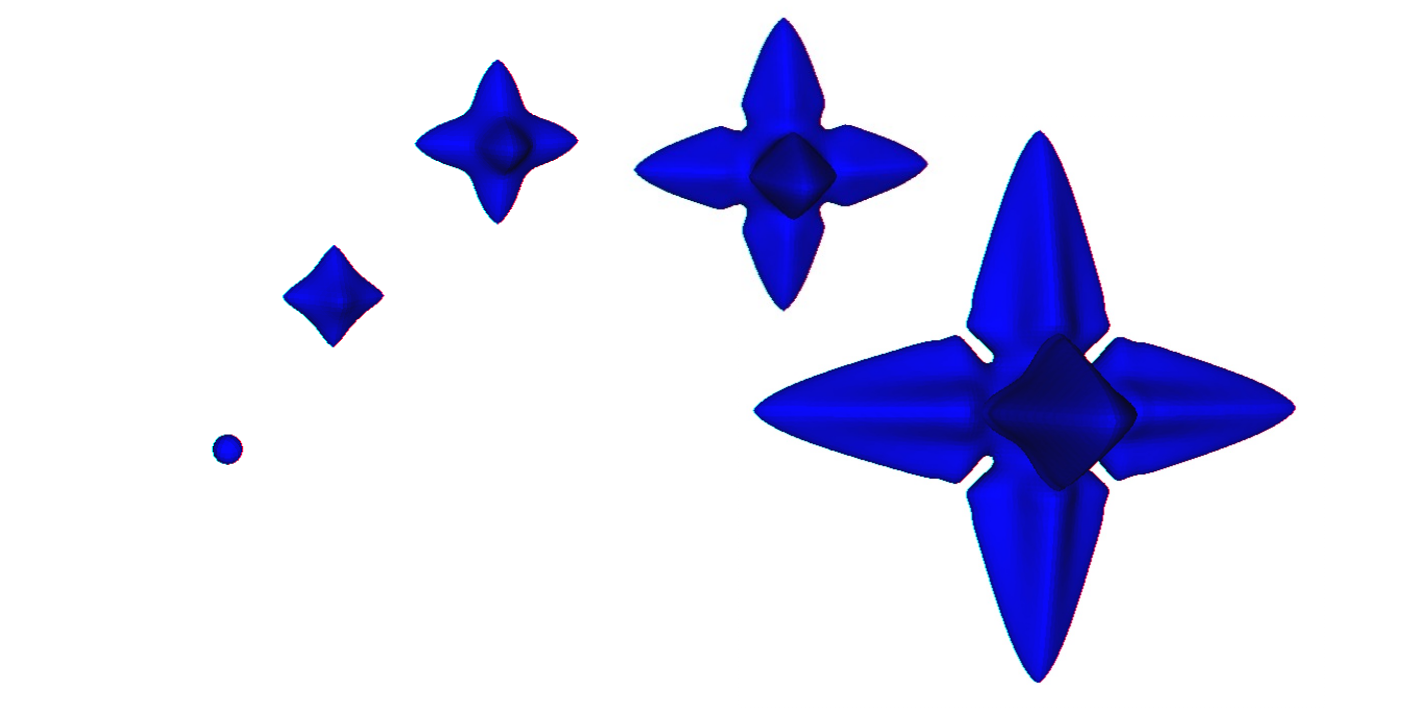}
    \caption{Time evolution of the free boundary with $\bar{\varepsilon}=0.001$ and $A=-3.0$ at $t=0$, $0.01$, $0.02$, $0.04$, and $0.08$.}
    \label{fig:stefan1}
\end{figure}

\begin{figure}[htbp]
    \centering
    \includegraphics[width=0.4\textwidth]{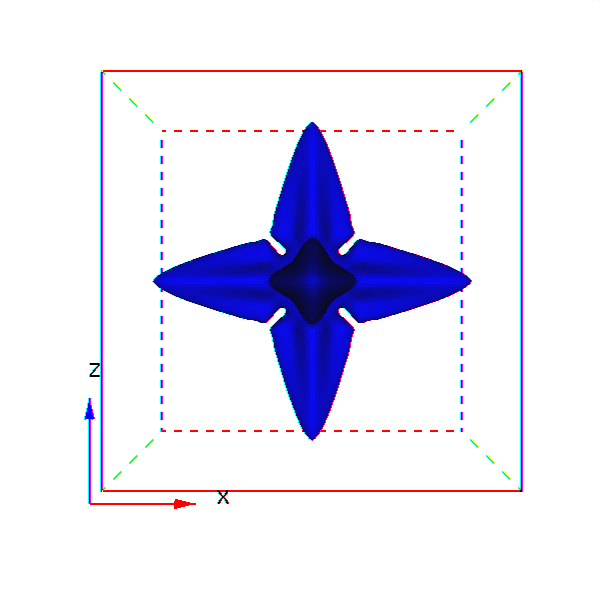}
    \includegraphics[width=0.4\textwidth]{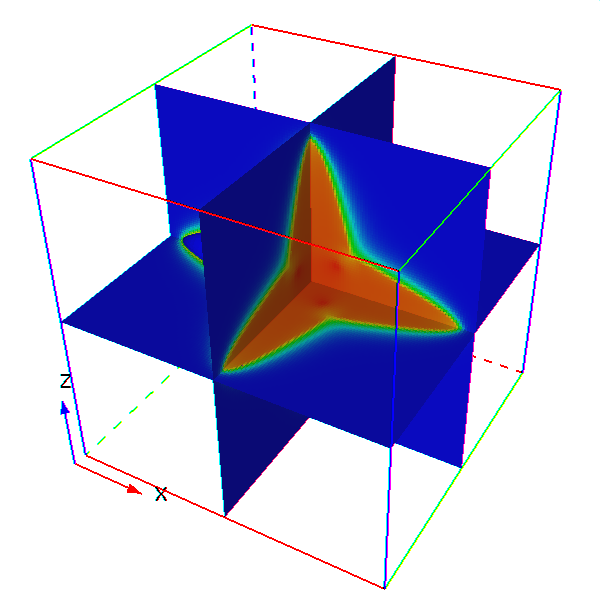}
    \caption{Front view of the free boundary (left) and slices of the temperature field (right).}
    \label{fig:stefan2}
\end{figure}

First, we choose $\bar{\varepsilon}=0.001$ and $A=-3.0$, so that dendrites are expected to grow in the six coordinate-axis directions.
In \Cref{fig:stefan1}, snapshots of the free boundary at $t=0, 0.01, 0.02, 0.04$, and $0.08$ are presented.
The morphology of the free boundary in the front view and the temperature field are also shown in \Cref{fig:stefan2}.
One can also observe that the symmetry of the free boundary is well-preserved by the proposed method.

\begin{figure}[htbp]
    \centering
    \includegraphics[width=0.3\textwidth]{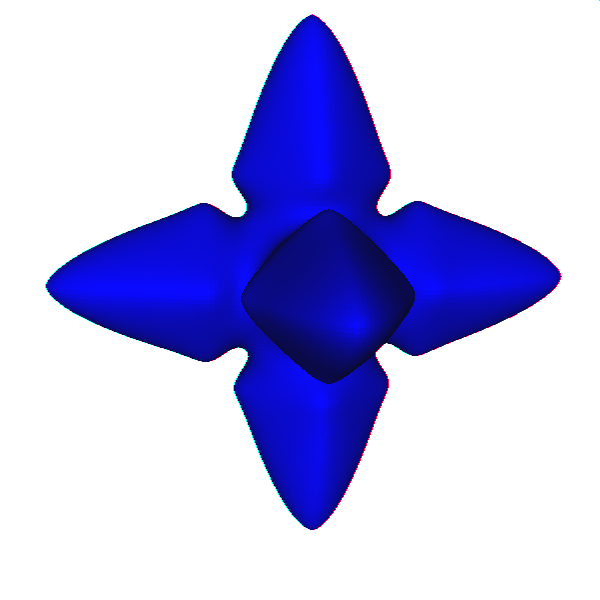}
    \includegraphics[width=0.3\textwidth]{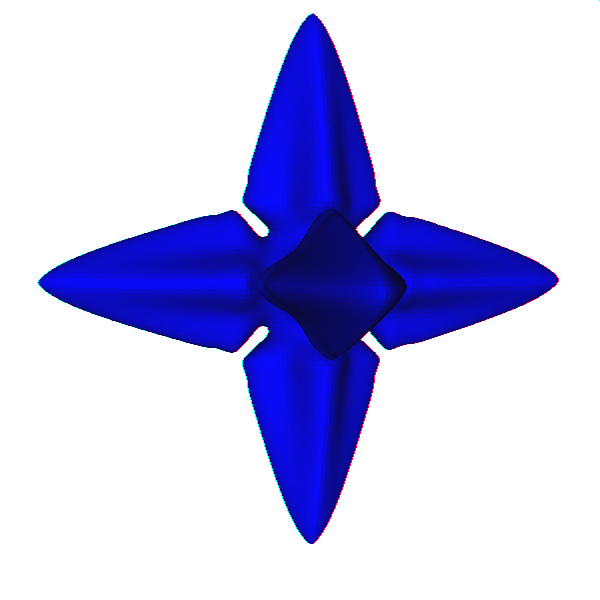}
    \includegraphics[width=0.3\textwidth]{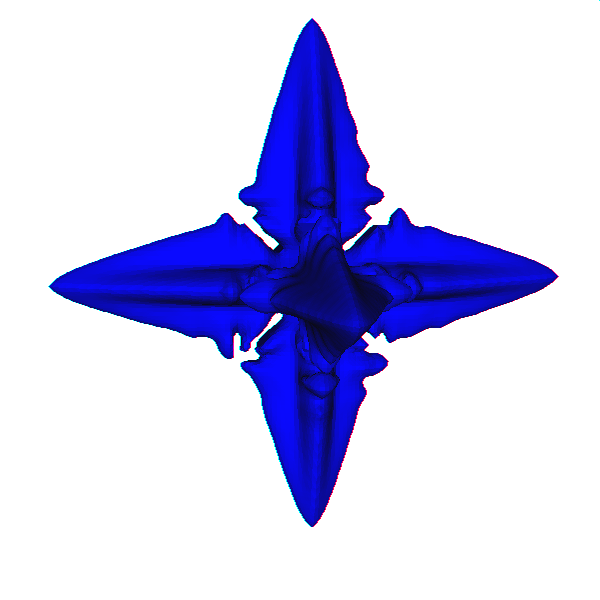}
    \caption{Morphologies of the free boundary using different parameters ($\bar{\varepsilon}=0.002$, $0.001$, and $0.0005$, from left to right).}
    \label{fig:stefan3}
\end{figure}

Then, different values of $\bar{\varepsilon}$ are chosen to vary the stabilizing effect on dendritic growth patterns.
Final morphologies of the free boundaries obtained with $\bar{\varepsilon} = 0.002$, $0.001$, and $0.0005$ are shown in \Cref{fig:stefan3}.
Larger values of $\bar{\varepsilon}$ lead to more stable results in the sense that they produce a smoother free boundary.

\begin{figure}[htbp]
    \centering
    \includegraphics[width=0.8\textwidth]{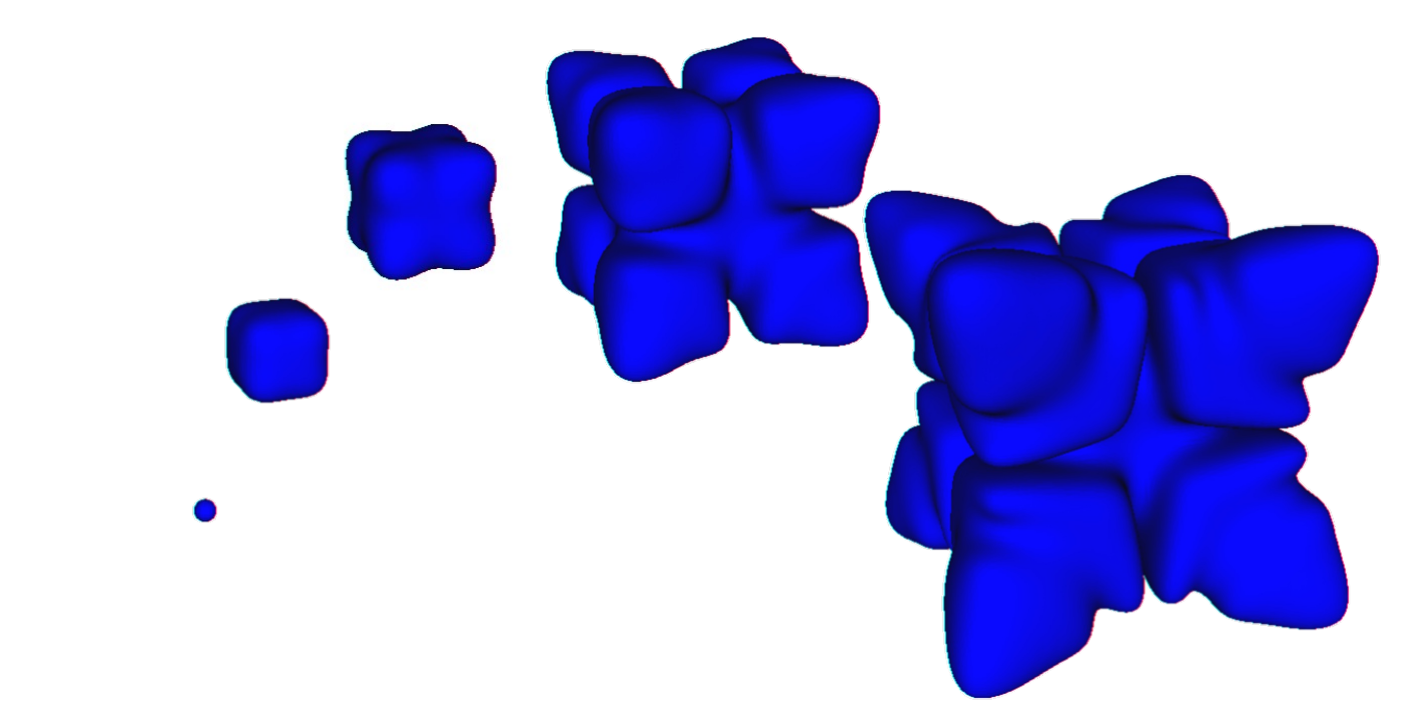}
    \caption{Time evolution of the free boundary with $\bar{\varepsilon}=0.005$ and $A=0.5$ at $t=0$, $0.05$, $0.1$, $0.3$, and $0.5$.}
    \label{fig:stefan4}
\end{figure}

We also choose $\bar{\varepsilon} = 0.005$ and $A=0.5$, so that the preferred growth directions are expected to be the eight vectors $(\pm 1, \pm 1, \pm1)$.
Snapshots of the free boundary at $t=0, 0.05, 0.10, 0.20$, and $0.50$ are presented in \Cref{fig:stefan4}.

\section{Discussion}\label{sec:discu}
In this paper, we develop ADI schemes for the heat equation, the reaction-diffusion equation, and the Stefan problem on arbitrary three-dimensional domains.
For problems on fixed domains, the proposed KFBI-ADI method achieves second-order accuracy in both space and time, and the modified DG scheme is unconditionally stable.
For problems on time-dependent domains, the level set-ADI method captures the free boundary with a level set formulation and efficiently solves the Stefan problem, enabling the simulation of complex dendritic growth patterns.

The proposed method is efficient in the sense that its computational complexity is linear in the number of degrees of freedom.
Furthermore, due to the dimension-splitting strategy of the ADI schemes, the method can be parallelized in a straightforward way, which can significantly accelerate the computation.

Since higher-order methods can attain a prescribed accuracy with fewer degrees of freedom, the design of higher-order extensions is a natural direction for future work.

This paper is mainly concerned with parabolic PDEs, whose discretizations share common features because the principal operator is the Laplacian.
For other types of PDEs, such as Maxwell's equations, which are hyperbolic, developing efficient and geometrically flexible ADI schemes is also important and will be explored in future work.

\appendix
\section{Stability analysis of the modified Douglas-Gunn ADI scheme}\label{sec:analysis}
Let $\lambda= \tau/h^2$ where $\tau$ is the time step and $h$ is the mesh size.
Denote by $\delta_{x}^{2}/h^2$, $\delta_{y}^{2}/h^2$, and $\delta_{z}^{2}/h^2$ the central difference operators for approximating $\partial_{xx}$, $\partial_{yy}$, and $\partial_{zz}$, respectively.
Let $u$ be the numerical solution to the heat equation.
The modified DG-ADI method can be written as
\begin{subequations}
\begin{align}
(1-\frac{\lambda}{2} \delta_{x}^{2}) u^{*} & =(1+\frac{\lambda}{2} \delta_{x}^{2}) u^{n}+\frac{3 \lambda}{2}(\delta_{y}^{2}+\delta_{z}^{2}) u^{n}-\frac{\lambda}{2}(\delta_{y}^{2}+\delta_{z}^{2}) u^{n-1}, \label{eqn:A1a} \\
(1-\frac{\lambda}{2} \delta_{y}^{2}) u^{* *} & =u^{*}+\frac{\lambda}{2} \delta_{y}^{2} u^{n-1}-\lambda \delta_{y}^{2} u^{n}, \label{eqn:A1b}\\
(1-\frac{\lambda}{2} \delta_{z}^{2}) u^{n+1} & =u^{* *}+\frac{\lambda}{2} \delta_{z}^{2} u^{n-1}-\lambda \delta_{z}^{2} u^{n}. \label{eqn:A1c}
\end{align}
\end{subequations}
Applying $(1-\frac{\lambda}{2} \delta_{x}^{2})$ to both sides of \eqref{eqn:A1b} yields
\begin{equation}\label{eqn:A2}
(1-\dfrac{\lambda}{2} \delta_{x}^{2})(1-\dfrac{\lambda}{2} \delta_{y}^{2}) u^{* *}=(1-\dfrac{\lambda}{2} \delta_{x}^{2}) u^{*}+\dfrac{\lambda}{2}(1-\dfrac{\lambda}{2} \delta_{x}^{2}) \delta_{y}^{2} u^{n-1}-\lambda(1-\dfrac{\lambda}{2} \delta_{x}^{2}) \delta_{y}^{2} u^{n}.
\end{equation}
Substituting \eqref{eqn:A1a} into \eqref{eqn:A2} yields
\begin{equation}
\begin{aligned}
(1-\dfrac{\lambda}{2} \delta_{x}^{2})(1-\dfrac{\lambda}{2} \delta_{y}^{2}) u^{* *}
&= u^{n}+\dfrac{\lambda}{2} \delta_{x}^{2} u^{n}+\dfrac{\lambda}{2} \delta_{y}^{2} u^{n} +\dfrac{3 \lambda}{2} \delta_{z}^{2} u^{n} \\
&-\dfrac{\lambda}{2} \delta_{z}^{2} u^{n-1}-\dfrac{\lambda^{2}}{4} \delta_{x}^{2} \delta_{y}^{2} u^{n-1}+\dfrac{\lambda^{2}}{2} \delta_{x}^{2} \delta_{y}^{2} u^{n},
\end{aligned}
\end{equation}
Similarly, applying $(1-\frac{\lambda}{2} \delta_{x}^{2})(1-\frac{\lambda}{2} \delta_{y}^{2})$ to \eqref{eqn:A1c} yields the recurrence relation for $u^{n+1}$, $u^n$, and $u^{n-1}$,
\begin{equation}\label{eqn:A4}
\begin{aligned}
& \quad (1-\dfrac{\lambda}{2} \delta_{x}^{2})(1-\dfrac{\lambda}{2} \delta_{y}^{2})(1-\dfrac{\lambda}{2} \delta_{z}^{2}) u^{n+1} \\
& =(1-\dfrac{\lambda}{2} \delta_{x}^{2})(1-\dfrac{\lambda}{2} \delta_{y}^{2}) u^{* *}+(1-\dfrac{\lambda}{2} \delta_{x}^{2})(1-\dfrac{\lambda}{2} \delta_{y}^{2})(\dfrac{\lambda}{2} \delta_{z}^{2} u^{n-1}-\lambda \delta_{z}^{2} u^{n}) \\
& =u^{n}+\dfrac{\lambda}{2} \delta_{x}^{2} u^{n}+\dfrac{\lambda}{2} \delta_{y}^{2} u^{n}+\dfrac{\lambda}{2} \delta_{z}^{2} u^{n}+\dfrac{\lambda^{2}}{2} \delta_{x}^{2} \delta_{y}^{2} u^{n}+\dfrac{\lambda^{2}}{2} \delta_{x}^{2} \delta_{z}^{2} u^{n} \\
& +\dfrac{\lambda^{2}}{2} \delta_{y}^{2} \delta_{z}^{2} u^{n}-\dfrac{\lambda^{2}}{4} \delta_{x}^{2} \delta_{y}^{2} u^{n-1}-\dfrac{\lambda^{2}}{4} \delta_{x}^{2} \delta_{z}^{2} u^{n-1}-\dfrac{\lambda^{2}}{4} \delta_{y}^{2} \delta_{z}^{2} u^{n-1} \\
& +\dfrac{\lambda^{3}}{8} \delta_{x}^{2} \delta_{y}^{2} \delta_{z}^{2} u^{n-1}-\dfrac{\lambda^{3}}{4} \delta_{x}^{2} \delta_{y}^{2} \delta_{z}^{2} u^{n} \\
& =(1+\dfrac{\lambda}{2} \delta_{x}^{2})(1+\dfrac{\lambda}{2} \delta_{y}^{2})(1+\dfrac{\lambda}{2} \delta_{z}^{2}) u^{n}+\dfrac{\lambda^{2}}{4}(\delta_{x}^{2} \delta_{y}^{2}+\delta_{x}^{2} \delta_{z}^{2}+\delta_{y}^{2} \delta_{z}^{2}) u^{n} \\
& -\dfrac{3 \lambda^{3}}{8} \delta_{x}^{2} \delta_{y}^{2} \delta_{z}^{2} u^{n}-\dfrac{\lambda^{2}}{4}(\delta_{x}^{2} \delta_{y}^{2}+\delta_{x}^{2} \delta_{z}^{2}+\delta_{y}^{2} \delta_{z}^{2}) u^{n-1} +\dfrac{\lambda^{3}}{8} \delta_{x}^{2} \delta_{y}^{2} \delta_{z}^{2} u^{n-1}.
\end{aligned}
\end{equation}
Next, the standard discrete Fourier analysis can be applied to study the $L_2$ stability of the ADI scheme.
Let $v_{i, j, k}^{n+1}=u_{i, j, k}^{n}$ and $\bm{U}_{i, j, k}^{n}=[u_{i, j, k}^{n}, v_{i, j, k}^{n}]$. Then \eqref{eqn:A4} can be reformulated as
\begin{equation}\label{eqn:A5}
    \bm{A}\bm{U}^{n+1} = \bm{B}\bm{U}^{n}
\end{equation}
where
$$
\bm{A} = 
\begin{pmatrix}
(1-\dfrac{\lambda}{2} \delta_{x}^{2})(1-\dfrac{\lambda}{2} \delta_{y}^{2})(1-\dfrac{\lambda}{2} \delta_{z}^{2}) & 0 \\
0 & 1
\end{pmatrix},
\quad
\bm{B} = 
\begin{pmatrix}
b_1 & b_2 \\
0 & 1
\end{pmatrix},
$$
$$ b_1=(1+\frac{\lambda}{2} \delta_{x}^{2})(1+\frac{\lambda}{2} \delta_{y}^{2})(1+\frac{\lambda}{2} \delta_{z}^{2})+\frac{\lambda^{2}}{4}(\delta_{x}^{2} \delta_{y}^{2}+\delta_{x}^{2} \delta_{z}^{2}+\delta_{y}^{2} \delta_{z}^{2})-\frac{3 \lambda^{3}}{8} \delta_{x}^{2} \delta_{y}^{2} \delta_{z}^{2}, $$
$$ b_2=-\frac{\lambda^{2}}{4}(\delta_{x}^{2} \delta_{y}^{2}+\delta_{x}^{2} \delta_{z}^{2}+\delta_{y}^{2} \delta_{z}^{2})+\frac{\lambda^{3}}{8} \delta_{x}^{2} \delta_{y}^{2} \delta_{z}^{2}. $$

Let us consider the single mode $\bm{U}_{ijk}^n = e^{\sqrt{-1}(i k_1  h + j k_2  h + k k_3  h)}\bm{V}^n$.
Then \eqref{eqn:A5} becomes
\begin{equation}
    \hat{\bm{A}} \bm{V}^{n+1} = \hat{\bm{B}} \bm{V}^{n},
\end{equation}
where
$$
\hat{\bm{A}} = 
\begin{pmatrix}
(1+\lambda \mu_{1})(1+\lambda \mu_{2})(1+\lambda \mu_{3}) & 0 \\
0 & 1
\end{pmatrix},
\quad
\hat{\bm{B}} = 
\begin{pmatrix}
\hat b_1 & \hat b_2 \\
0 & 1
\end{pmatrix},
$$
$$ \hat b_1=(1-\lambda \mu_{1})(1-\lambda \mu_{2})(1-\lambda \mu_{3})+\lambda^{2}(\mu_{1} \mu_{2}+\mu_{1} \mu_{3}+\mu_{2} \mu_{3})+3 \lambda^{3} \mu_{1} \mu_{2} \mu_{3}, $$
$$ \hat b_2=-\lambda^{2}(\mu_{1} \mu_{2}+\mu_{1} \mu_{3}+\mu_{2} \mu_{3})-\lambda^{3} \mu_{1} \mu_{2} \mu_{3}, $$
$$ \mu_{1}=2 \sin ^{2} \frac{k_{1} h}{2}, \mu_{2}=2 \sin ^{2} \frac{k_{2} h}{2}, \mu_{3}=2 \sin ^{2} \frac{k_{3} h}{2}. $$
The characteristic equation of the amplification matrix $\hat{\bm{A}}^{-1} \hat{\bm{B}}$ is given by
\begin{equation}\label{eqn:char-eqn}
    \Lambda^{2} - b \Lambda - c = 0,
\end{equation}
where
$$ b  =\dfrac{(1-\lambda \mu_{1})(1-\lambda \mu_{2})(1-\lambda \mu_{3})+\lambda^{2}(\mu_{1} \mu_{2}+\mu_{1} \mu_{3}+\mu_{2} \mu_{3})+3 \lambda^{3} \mu_{1} \mu_{2} \mu_{3}}{(1+\lambda \mu_{1})(1+\lambda \mu_{2})(1+\lambda \mu_{3})},$$
$$ c  =\dfrac{-\lambda^{2}(\mu_{1} \mu_{2}+\mu_{1} \mu_{3}+\mu_{2} \mu_{3})-\lambda^{3} \mu_{1} \mu_{2} \mu_{3}}{(1+\lambda \mu_{1})(1+\lambda \mu_{2})(1+\lambda \mu_{3})}.$$
Recall that the moduli of the roots of \eqref{eqn:char-eqn} are no greater than $1$ if and only if $|b|\leq 1-c$ and $|c|\leq 1$.
It is obvious that the inequality $|c| \leq 1$ holds and we only need to verify $c-1\leq b \leq 1-c$.
Since we have
\begin{equation}
\begin{aligned}
 &-\lambda^{2}(\mu_{1} \mu_{2}+\mu_{1} \mu_{3}+\mu_{2} \mu_{3})-\lambda^{3} \mu_{1} \mu_{2} \mu_{3} \\
& \leq (1+\lambda \mu_{1})(1+\lambda \mu_{2})(1+\lambda \mu_{3})+(1-\lambda \mu_{1})(1-\lambda \mu_{2})(1-\lambda \mu_{3}) \\
& + \lambda^{2}(\mu_{1} \mu_{2}+\mu_{1} \mu_{3}+\mu_{2} \mu_{3})+3 \lambda^{3} \mu_{1} \mu_{2} \mu_{3} \\
& = 2+3 \lambda^{2}(\mu_{1} \mu_{2}+\mu_{1} \mu_{3}+\mu_{2} \mu_{3})+3 \lambda^{3} \mu_{1} \mu_{2} \mu_{3},
\end{aligned}
\end{equation}
then it gives $c-1\leq b$.
The inequality from the other side can be seen from
\begin{align*}
&(1-\lambda \mu_{1})(1-\lambda \mu_{2})(1-\lambda \mu_{3})+\lambda^{2}(\mu_{1} \mu_{2}+\mu_{1} \mu_{3}+\mu_{2} \mu_{3})+3 \lambda^{3} \mu_{1} \mu_{2} \mu_{3} \\
&\leq(1+\lambda \mu_{1})(1+\lambda \mu_{2})(1+\lambda \mu_{3}) +\lambda^{2}(\mu_{1} \mu_{2}+\mu_{1} \mu_{3}+\mu_{2} \mu_{3}) +\lambda^{3} \mu_{1} \mu_{2} \mu_{3},
\end{align*}
which is equivalent to $0 \leq \lambda(\mu_{1}+\mu_{2}+\mu_{3})$.
Therefore, the spectral radius of the amplification matrix $\hat{\bm{A}}^{-1} \hat{\bm{B}}$ is no greater than $1$, and the scheme is unconditionally stable.

\section*{Acknowledgement}
W.~Ying is supported by the National Natural Science Foundation of China in the Division of Mathematical Sciences (Project No.~12471342) and the fundamental research funds for the central universities of China.

\bibliographystyle{elsarticle-num}
\bibliography{references}

\end{document}